%% file: main_arxiv.tex
\pdfoutput=1
\documentclass{article}
\usepackage{amsmath,amssymb,dsfont}
\usepackage[margin=1in]{geometry}
\usepackage{cite}
\usepackage{xcolor}
\usepackage{amsfonts}
\usepackage{graphicx}
\usepackage{tikz}
\usepackage{lipsum}
\usepackage{MnSymbol}
\usepackage{mathtools}
\usepackage{algorithm}
\usepackage{algpseudocode}
\usepackage{amsthm}
\usepackage[mathscr]{euscript}
\def\spoke#1#2{
\begin{scope}[shift={#1}, rotate=#2]
    \draw[thick] (0,-.15) -- (0,-1);
\end{scope}
}

\providecommand{\red}[1]{\textcolor{black}{#1}}

\providecommand{\R}{\ensuremath \mathbb{R}}
\providecommand{\N}{\ensuremath \mathbb{N}}

\providecommand{\ip}[1]{\ensuremath \left\langle #1\right\rangle}
\providecommand{\spt}{\textrm{spt}}

\newtheorem{thm}{Theorem}

\newtheorem{lem}[thm]{Lemma}

\newtheorem{assum}[thm]{Assumption}

\theoremstyle{definition}
\newtheorem{defn}[thm]{Definition}
\newtheorem{exmp}[thm]{Example}
\newtheorem{rem}[thm]{Remark}

\newenvironment{pf}{\begin{proof}}{\end{proof}}

\newcommand*\circled[1]{\tikz[baseline=(char.base)]{
            \node[shape=circle,draw,inner sep=.pt] (char) {#1};}}

\begin{document}
\title{Convex Computation of the Reachable Set for Hybrid Systems with Parametric Uncertainty}
\author{Shankar Mohan\thanks{Shankar Mohan (\texttt{elemsn@umich.edu}) is with the Department of Electrical Engineering and Computer Science, University of Michigan, Ann Arbor, MI 48109}, Victor Shia\thanks{Victor Shia (\texttt{vshia@eecs.berkeley.edu}) is with the Department of Electrical Engineering and Computer Sciences, University of California at Berkeley, Berkeley, CA 94720}, Ram Vasudevan\thanks{Ram Vasudevan (\texttt{ramv@umich.edu}) is with the Department of Mechanical Engineering, University of Michigan, Ann Arbor, MI 48109}}
\maketitle
  \begin{abstract}
  	To verify the correct operation of systems, engineers need to determine the set of configurations of a dynamical model that are able to safely reach a specified configuration under a control law. Unfortunately, constructing models for systems interacting in highly dynamic environments is difficult. This paper addresses this challenge by presenting a convex optimization method to efficiently compute the set of configurations of a polynomial hybrid dynamical system that are able to safely reach a user defined target set despite parametric uncertainty in the model. This class of models describes, for example, legged robots moving over uncertain terrains.	The presented approach utilizes the notion of occupation measures to describe the evolution of trajectories of a nonlinear hybrid dynamical system with parametric uncertainty as a linear equation over measures whose supports coincide with the trajectories under investigation. This linear equation with user defined support constraints is approximated with vanishing conservatism using a hierarchy of semidefinite programs \red{that are each} proven to compute an \red{inner/}outer approximation to the set of initial conditions that can reach the user defined target set safely in spite of uncertainty. The efficacy of this method is illustrated on a collection of \red{six} representative examples.
  \end{abstract}

  \input{sections/sec_intro}
  \input{sections/sec_notations}
  \input{sections/sec_PF}
  \input{sections/sec_relax}
  \input{sections/sec_extensions}
  \input{sections/sec_examples}
  \input{sections/sec_conc}
  \appendix
  \input{sections/sec_app}

\end{document}

%% file: sections/sec_intro.tex
\section{Introduction}
Computing the set of configurations that are able to safely reach a desired configuration is critical to ensuring the correct performance of a system in dynamic environments where deviations from planned behavior are to be expected.
Many methods have been proposed to efficiently compute this set that is generally referred to as the \emph{backwards reachable set} for deterministic systems.
Unfortunately, the effect of intermittent contact with the world, especially in fluctuating environments, is demanding to model deterministically.
A roboticist, for example, may be tasked with ensuring that a control for a legged robot beginning from an initial configuration is able to safely reach a desired goal; however, limitations in sensing or environment variability may render exact modeling of terrain height or friction impossible.
The development of numerical tools to tractably compute the backwards reachable set of dynamical systems undergoing contact, or \emph{hybrid dynamical systems}, with parametric uncertainty while providing systematic guarantees has been challenging due to the difficulty of efficiently accounting for the uncertainty.
\par
Given the potential utility of the set of configurations that are able to reach a target set despite parametric uncertainty, called the \emph{uncertain backwards reachable set}, many researchers have attempted to develop numerical tools to approximate this set. For instance, in \cite{Chesi2004}, given a Lyapunov function, the largest ellipsoidal level-set contained in the {\em uncertain backwards reachable set} (where the terminal set is the origin) for polynomial systems is computed by formulating an eigenvalue problem. While the method developed is able to address systems with either time-varying or constant uncertainties, the computed ellipsoid is dependent on the provided Lyapunov function and hence is not guaranteed to be the {\em uncertain backwards reachable set}.
\par
Several other researchers have attempted to utilize this set to construct controllers for legged robots that are able to walk over terrains of varying heights \cite{byl2008metastable,dai2012optimizing,griffin2015,saglam2013switching}.
Their approach has relied on discretizing the height of the terrain or selecting specific terrain profiles while constructing a safe controller across these specified heights, which verifies the performance of the controller only at those specific heights.
Moreover, these approaches are unable to account for uncertainty associated with imperfect knowledge of terrain friction or parameters affecting the continuous dynamics.
\par
Other researchers have developed tools to outer approximate the {\em uncertain backwards reachable} for linear systems with uncertain parameters using a variety of approaches \cite{girard2005reachability,althoff2008reachability}. These methods can be extended to nonlinear hybrid systems, but can require the introduction of a large number of discrete states to represent the nonlinear behavior or require overly conservative estimates of potential uncertainty; this approach of using zonotopes to generate an envelope for system trajectories has been more recently explored in \cite{Maiga2015}. More generally, Hamilton-Jacobi Bellman based approaches have also been applied to compute the {\em uncertain backwards reachable set} for nonlinear systems with arbitrary uncertainty affecting the state at any instance in time \cite{tomlin2003computational}.
These approaches solve a more general problem, but rely on state space discretization which can be prohibitive for systems of dimension greater than four without relying upon specific system structure \cite{maidens2013lagrangian}.
\par
This paper leverages a method developed in several recent papers \cite{henrion2014convex,majumdar2014convex,shia2014convex,Sloth2015} that describe the evolution of trajectories of a deterministic hybrid dynamical system using measures, to describe the evolution of a hybrid dynamical system with parametric uncertainty as a linear equation over measures. As a result of this characterization the {\em uncertain backwards reachable set} can be computed as the solution to an infinite dimensional linear program over the space of nonnegative measures. To compute an approximate solution to this infinite dimensional linear program, a sequence of finite dimensional relaxed semi-definite programs are constructed that satisfy an important property:
each solution to this sequence of semi-definite programs is an outer approximation to the {\em uncertain backwards reachable set} with asymptotically vanishing conservatism.
The approach is most comparable to those that check Lyapunov's criteria for stability via sums-of-squares programming to verify the safety of a system \cite{prajna2004safety,Posa2015}.
In contrast to these approaches, the algorithm described in this paper does not require solving a bilinear optimization problem that requires feasible initialization and allows for more general descriptions of the parametric uncertainty in the model.
\par
The remainder of the paper is organized as follows:
Section~\ref{sec:preliminaries} introduces the notation used in the remainder of the paper, the class of systems under consideration, and the backwards reachable set problem under parametric uncertainty;
Section~\ref{sec:prob} describes how the backwards reachable set under parametric uncertainty is the solution to an infinite dimensional linear program;
Section~\ref{sec:implementation} constructs a sequence of finite dimensional semidefinite programs that outer approximate the infinite dimensional linear program with vanishing conservatism;
Section~\ref{sec:ext} provides extensions to the contents in Section~\ref{sec:prob} by presenting amongst others, a converging inner approximations of the BRS;
Section~\ref{sec:examples} describes the performance of the approach on a set of examples;
and, Section~\ref{sec:conclusion} concludes the paper.

%% file: sections/sec_notations.tex
  \section{Preliminaries}
\label{sec:preliminaries}
  This section defines the notation, the class of systems, and the problem considered throughout this paper. The reader is directed to \cite{folland2013real,Bogachev2007,Bogachev2007a} for an introduction to some of the measure theoretic concepts utilized in this paper.
  \subsection{Notation}
  In the remainder of this text the following notation is adopted: sets are italicized and capitalized (ex. $K$).
The boundary of a set $K$ is denoted by $\partial K$.
Finite truncations of the set of natural numbers are expressed as \mbox{$\N_n:=\{1,\ldots,n\}$}.
The set of continuous functions on a compact set $K$ are denoted by $\mathcal C(K)$.
The ring of polynomials in $x$ is denoted by $\R[x]$, and the degree of a polynomial is equal to degree its largest multinomial; the degree of the multinomial $x^\alpha,\,\alpha\in \N^n$ is $|\alpha|=\|\alpha\|_1$; and $\R_d[x]$ is the set of polynomials in $x$ with maximum degree $d$.

The dual to $\mathcal C(K)$ is the set of Radon measures on $K$, denoted as $\mathcal M(K)$, and the pairing of $\mu\in \mathcal M(K)$ and $v\in \mathcal C(K)$ is:
  \begin{align}
  \ip{\mu,v}=\int_{K}v(x)\,d\mu(x).
  \end{align}
We denote the nonnegative Radon measures by ${\mathcal M}_+(K)$.
The space of Radon probability measures on $K$ is denoted by ${\mathcal P}(K)$.
The Lebesgue measure is denoted by $\lambda$.
Finally, supports of measures, $\mu$, are identified as $\spt(\mu)$.

\subsection{Quasi-Uncertain Hybrid Systems}
\label{ssec:hybrid}
We next define the class of uncertain hybrid systems considered throughout the remainder of the paper; the definition is an adaptation of the description in \cite{Burden2015}.
\begin{defn}
\label{def:system}
  A `quasi-uncertain' hybrid system is a tuple \mbox{$\mathcal H=(\mathcal J,\mathcal E,\mathcal D,\red{\Gamma},\mathcal F,\mathcal G,\mathcal R)$}, where
  \begin{itemize}
    \item $\mathcal J$ is a finite set of indices of discrete states in of $\mathcal H$ (discrete states will also be referred to as modes);
    \item $\mathcal E\subset \mathcal J\times \mathcal J$ is a set of two-tuples describing directed edges;
    \item $\mathcal D = \{D_j\}_{j \in {\mathcal J}}$ is the set of domains where each $D_j$ is a compact $n_j$-dimensional manifold with boundary where $n_j \in \N$;
    \item $\Gamma = \{\mu_{\theta_j} \}_{j \in {\mathcal J}}$ where $\mu_{\theta_j}\in \mathcal P(\Theta_j)$ describes the uncertainty associated with discrete state $j \in {\mathcal J}$ with $\Theta_j$ being a compact set;
    \item $\mathcal F=\{\tilde f_j\}_{j\in \mathcal J}$ where $\tilde{f}_j: D_j \times \Theta_j \to D_j$ is a Lipschitz continuous function \red{(in all variables)} describing the dynamics on $D_j$;
    \item \red{$\mathcal G=\{G_e\}_{e \in {\mathcal E}}$ is the set of guards where each $G_{(j,j')} \subset \partial D_j$ is a guard in domain $j \in {\mathcal J}$ that defines a transition from mode $j$ to mode $j' \in {\mathcal J}$};
    \item $\mathcal R = \{R_e\}_{e \in {\mathcal E}}$ is the set of continuous reset maps, where each map is a continuously differentiable injection $R_{(j,j')}: G_{(j,j')} \to D_{j'}$.
  \end{itemize}
\end{defn}
The definition of reset maps needs some clarification. Let $A\subset G_{(j,k)}$ be a set such that $R_{(j,k)}(A)\subset G_{(k,l)}$. In this case, we will not be altering the system behaviour by identifying a new guard (or coalescing with an existing guard)  $G_{(j,l)}$ such that $A\subset G_{(j,l)}$ and re-defining $G_{(j,k)}:=G_{(j,k)}\backslash A$. Hence, without loss of generality, we require that
\begin{align}
    R_e(G_e) \bigcap\limits_{ \shortstack{$a\in \mathcal E$\\$a\ne e$} } G_a =\emptyset,\phantom{8}\forall e\in \mathcal E.
    \label{eq:reset}
\end{align}
\par
To avoid any ambiguity during transitions between discrete states, we assume the following:
\begin{assum}
    In each discrete state, the guards are mutually exclusive; i.e.
    \begin{align}
     G_{(i,j)}\cap G_{(i,k)}=\emptyset,\phantom{8}\forall (i,j),(i,k)\in \mathcal E, \forall j\ne k
    \end{align}
\end{assum}
In addition, the systems are not allowed to undergo infinite mode transitions in any finite time-interval.
\begin{assum}
  $\mathcal H$ has no zeno execution.
  \label{assump:zeno}
\end{assum}
Algorithm~\ref{alg:execution} describes the finite-time execution, $[0,T]$, of a hybrid system, $\mathcal H$, as in Defn.~\ref{def:system} as follows:
Suppose that the system enters mode $j$ at time $t$ at location $x \in D_j$.
Recall that the dynamics in this domain, $\tilde f_j$, are a function of a random parameter drawn from the distribution $\mu_{\theta_j}$; let this random variable take the value $\theta$.
The trajectory of the hybrid system beginning at time $t$ at $x$ is then given by any absolutely continuous function that satisfies the differential equation $\tilde f_j$ with a fixed $\theta$ as described in Steps~\circled{5}\&\circled{6}.
This trajectory evolves until either the time evolution passes $T$ or the trajectory arrives at a guard, whichever happens first.
Steps~\circled{7}--\circled{11} isolates the first hitting-time of a guard of mode $j$ and resets \red{the solution trajectory} to a new mode whereafter the same procedure is repeated until $t=T$.

\begin{algorithm}[!t]
\small
\begin{algorithmic}[1]
 \State {\bf Initialization:} $t=0,\,j\in \mathcal J,\,x_0\in D_j,\,x(0)=x_0$\;
 \State{\bf While} 1 {\bf do}\;
 \State\hspace{.2in}{\em Let} $\theta$ be drawn according to $\mu_{\theta_j}$\;
 \State\hspace{.2in}{\em Let} $\gamma\colon [t,T]\rightarrow D_j$, absolutely continuous st.\;
 \State\hspace{.4in}$\dot \gamma(s)=\tilde f(\gamma(s),\theta)$ $\lambda_t^{\text{\tiny a}}$-a.e., $s\in [t,T]$\;
 \State\hspace{.4in}$\gamma(t)=x(t)$\;
 \State\hspace{.2in}$\Lambda_{(j,t)}:=\{r\in [t,T]| \exists (j,k)\in \mathcal E \text{ st. } (\gamma(r),\theta)\in G_{(j,k)}\}$\;
 \State\hspace{.2in}{\bf If} {$\Lambda_{(j,t)}\ne \emptyset$} {\bf then}\;
 \State\hspace{.4in}$t':=\min \Lambda_{(j,t)}$, $k$ st. $\gamma(t')\in G_{(j,k)}$\;
 \State\hspace{.4in}$x(s)\leftarrow \gamma(s)$, $\forall s\in [t,t')$\;
 \State\hspace{.4in}$t\leftarrow t',\,x(t')\leftarrow R_{(j,k)}(\gamma(t')),\,j\leftarrow k$
 \State\hspace{.2in}{\bf else}\;
 \State\hspace{.4in}$x(s)=\gamma(s),\,\forall s\in [t,T]$\;
 \State\hspace{.4in}Stop\;
 \State\hspace{.2in}{\bf end}\;
 \State{\bf end}\;
 \end{algorithmic}
\caption{Execution of $\mathcal H$}
 \label{alg:execution}
 $^a$where $\lambda_t$ is the Lebesgue measure on $[t,T]$
\end{algorithm}
Note that the uncertainty does not evolve with time; the uncertainty only changes when the \red{state trajectory is reset.}
This class of systems is rich as is illustrated by the following two examples: a simple 1D pedagogical example and a 2D representative of walking models.
\begin{exmp}[1-D Quasi-Uncertain Linear System]
\label{example:1D}
Consider a quasi-uncertain linear system evolving with dynamics:
\begin{align}
	\dot x = -0.7x+0.2\theta-0.1,\phantom{8} \forall x\in D_1
\end{align}
where \red{$\theta$ is an unknown parameter affecting the dynamics. The system state evolves on $D_1 = [-1,1]$ and the value of $\theta$ is drawn from a uniform distribution on $\Theta_1=[0.2,1]$. The uncertain parameter can be thought of as having arisen due to structural modeling errors, or as a result of reducing a system with time-scale separation.}
\end{exmp}
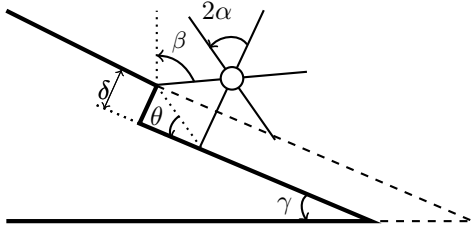
\begin{figure}[!t]
\centering
  \begin{tikzpicture}[scale=1]
    \draw[-,ultra thick] (0,2) to (2,1);
    \draw[dashed, thick] (2,1) -- (6.2,-0.8)-- (4.85,-.8);
    \draw[ultra thick] (2,1) to (1.77,0.5) -- (4.85,-0.8)--(0,-0.8);
    \draw[thick] (3,1.1) circle (.15cm);
    \spoke{(3.0,1.1)}{-85};
    \spoke{(3.0,1.1)}{-25};
    \spoke{(3.0,1.1)}{35};
    \spoke{(3,1.1)}{95};
    \spoke{(3,1.1)}{155};
    \spoke{(3,1.1)}{-145};
    \draw[dotted,thick] (2,1) -- (2,2);
    \draw[dotted, thick] (2,1) -- (2.63,0.12);
    \draw[->,thick] (2.3, 0.6) to [out=-145, in =115] (2.2,0.3) ;
    \node at (2.0,0.6) {$\theta$};
    \draw[->,thick] (4,-.45) to  [out=-145, in =125] (4,-0.8);
    \node at (3.7,-0.6) {$\gamma$};
    \draw[->, thick] (3.2,1.55) to [out=125, in =55] (2.7,1.5);
    \node at (2.8,2) {$2\alpha$};
    \draw[->,thick] (2.5,1.05) to [out=115, in =0] node[pos=0.5, above]{\small $\beta$} (2,1.4) ;
    \draw[dotted,thick] (1.77,0.5) -- (1.15,.765);
    \draw[<->] (1.3,.7)  -- node [pos=0,above] {$\delta$} (1.53,1.21);
    \draw[dotted,thick] (1.77,0.5) -- (1.15,.765);
    \draw[<->] (1.3,.7)  -- node [pos=0,above] {$\delta$} (1.53,1.21);
  \end{tikzpicture}
  \caption{Schematic of the rimless wheel with $\theta$ describing the effect of unknown terrain.}
  \label{fig:rw_schematic}
\end{figure}

\begin{exmp}[\red{Rimless Wheel on an uneven terrain}]
\label{example:rw}
The rimless wheel---constituted by a massless axle to which $n$ (angularly) equidistant spokes are connected---is a simple model of legged locomotion \cite{Posa2015,shia2014convex}. Figure~\ref{fig:rw_schematic} presents a schematic of a rimless wheel---with spokes separated by an angle $2\alpha$---rolling down an inclined plane. The rimless wheel is a hybrid system consisting of one mode; every time the spoke makes contact with the surface of the inclined plane, the system undergoes a reset. The continuous dynamics of the rimless wheel are:
\begin{align}
    \ddot\beta= \sin(\beta)
\end{align}
where $\beta$ is the angle between the vertical (which is defined as the line that is perpendicular to the base of the inclined plane) and the pivoting spoke.
Once the swinging spoke makes contact with the terrain, the states are reset as:
\begin{align}
R_{(1,1)}(\beta^-,\dot \beta^-)=\begin{bmatrix}
    2\gamma-\beta^-&
    \cos(2\alpha)\,\dot\beta^-
  \end{bmatrix}^T.
\end{align}
For a rimless wheel rolling down an inclined plane with flat terrain, at the instance when the swinging spoke makes contact with the ground, $\beta = \gamma + \alpha$. To encode the uncertainty due to terrain height, suppose the rimless wheel encounters a step of size $\delta$, then if we let $\theta =  \arcsin \left( \frac{\delta}{2 l \sin \alpha} \right)$, the guard is defined as:
\begin{align}
G_{(1,1)}=\{ (\beta, \dot{\beta},\red{\theta})\mid \beta=\gamma+\alpha+\theta\}.
\end{align}
Observe that as the rimless wheel continues to roll, the terrain is allowed to change since the random variable $\theta$ is allowed to take a distinct value after each contact with the ground.
\end{exmp}
\subsection{Problem Description}
{The objective of this work is to estimate the largest set of initial conditions from which all state trajectories of $\mathcal H$, regardless of any encountered uncertainty, reach a terminal set by a pre-specified time, $T$.}
To formalize the definition of this \emph{uncertain backwards reachable set}, we denote the terminal set as $X_T$ and its projection into each mode by $X_{(T,j)}$, which we assume is compact.
For convenience, we \red{let} $\mathcal T=[0,\,T]$.
We define the uncertain backwards reachable set mode-wise:
\begin{align}
\begin{aligned}
     X_{(0,j)}=\{x_0\in D_j \mid\ &\forall x: [0,T] \xrightarrow{\text{Alg.~\ref{alg:execution}}} {\mathcal D},& \text{with } x(0) = x_0, \, x(T) \in X_T \}
     \label{eq:brs}
\end{aligned}
\end{align}
The uncertain backwards reachable set is then defined as $X_0 = \{X_{(0,j)} \}_{j\in {\mathcal J}}$.
Observe that by definition, all initial conditions originating in any member of $X_0$ must reach $X_T$ at time $T$ regardless of mode transitions and uncertainty encountered along the way.
\par
\red{Lastly, we assume that mode transitions take place instantaneously. As a consequence, $X_T$ and guards must have empty intersections.
\begin{assum}
\label{assum:overlap}
  The guards and terminal set are mutually exclusive.
\end{assum}
}
\begin{rem}
\label{rem:guard_boundary}
  In describing the execution of the hybrid system, $\mathcal H$, it is assumed that guard is a subset of $\partial D$. There is no loss of generality in making this assumption, since it is always possible to define new modes (by partitioning existing modes) such that the guard lies at the boundary.
\end{rem}

%% file: sections/sec_PF.tex
\section{Problem Formulation}
\label{sec:prob}
In this section, we present a pair of dual infinite dimensional linear programs that compute the uncertain backwards reachable set.
Critically, note that despite the uncertainty being drawn from a distribution at the arrival into each mode, it remains constant throughout that mode.
As a result, this unknown parameter can be appended to the dynamics of every mode $j$ and treated as a portion of the state-space:
\begin{align}
f_j=\begin{bmatrix}
  \tilde f_j^T&\mathbf{0}^T_{n_{\theta_j}}
\end{bmatrix}^T.
\label{eq:ode}
\end{align}
To address the problem of estimating the {\em uncertain backwards reachable set} (BRS), we rely on the notion of \emph{occupation measures}, first introduced in \cite{Pitman1977}, to transform the hybrid nonlinear dynamics of the system into a set linear dynamics over measures that can more readily be solved.
For instance, suppose the system enters mode $j$ at $\tau_k$ with the states being initialized as $x(\tau_k)=x_0$ and $\theta(\tau_k)=\theta$,
the occupation measure, \mbox{$\mu_j(\cdot\mid \tau_k,x_0,\,\theta)\in \mathcal M_+(\mathcal T\times D_j\times \Theta_j)$}, is defined as:
\begin{align}
\hspace*{-1mm}\mu_j(A\times B\times C| \tau_k ,x_0,\theta)=\hspace*{-1.25mm}\int\limits_{0}^T \hspace*{-1.25mm} I_{A\times B\times C}(t,x(t|\tau_k,x_0,\theta),\theta)dt.
\end{align}
Note that the following relation between the Lebesgue measure on $\mathcal T$\red{, $\lambda$}, and $\mu_j(\cdot\mid \tau_k,x_0,\theta)$ holds for all $v \in C(\mathcal T\times D_j\times \Theta_j)$:
\begin{align}
\ip{\mu_j(\cdot\mid \tau_k,x_0,\theta),v}=\ip{\lambda,v(t,x(t\mid \tau_k,x_0,\theta),\theta)}.
\label{eq:mu_lambda}
\end{align}
The \emph{occupation measure}, as defined, is a conditional measure -- conditioned on the arrival-time and initial values of the states in that mode.
To consider a set of possible arrival-times and initial conditions, we define the \emph{average occupation measure} by integrating the conditional occupation measure against a measure on the set of possible initial conditions of the mode, $\mu_{s_j} \in  {\mathcal M}_+(\mathcal T\times D_j\times \Theta_j)$:
\begin{align}
\mu_j(A\times B\times C)= \hspace*{-5mm}\int\limits_{\mathcal T\times D_j\times \Theta_j}\hspace*{-5mm}\mu_j(A\times B\times C\mid \tau_k,x_0,\theta)\,d \mu_{s_j}.
\label{eq:mu_avg}
\end{align}
\red{The \emph{occupation measure} is to be interpreted as the total time that a particular solution trajectory (with a given initial conditions of states, uncertainty and initial time) spends in a set of interest. In contrast, the \emph{average occupation measure} is the total time spent by all solution trajectories whose initial conditions are in a given set, in a set of interest.}
\par
Observe that by definition, the uncertain variables are independent of the states' initial conditions; hence $\mu_{s_i}\in \mathcal M_+(\mathcal T\times D_j\times \Theta_j)$ is expressible as a product measure:
\begin{align}
\mu_{s_j}=\bar\mu_{0_j}\otimes \mu_{\theta_j},\label{eq:mu_sj}
\end{align}
where $\bar \mu_{0_j}\in \mathcal M_+(\mathcal T\times D_j)$ is a measure describing the set of initial conditions of the states, and $\mu_{\theta_j}\in \mathcal M_+(\Theta_j)$ is as in the definition of $\mathcal H$. Lemma~\ref{lemma:initial_measure} completes the characterization of the \emph{initial measure}, $\mu_{s_j}$.
\par
Similarly, measures on terminals sets, $\mu_{T_j}\in \mathcal M_+(X_{(T,j)}\times \Theta_j)$:
\begin{align}
\mu_{T_j}(A\times B)= \hspace*{-4mm}\int\limits_{\mathcal T\times D_j\times \Theta_j}\hspace*{-4mm}I_{A\times B}(x(T\mid \tau_k,x_0,\theta),\theta)\,d\mu_{s_{j}},
\end{align}
and guards, $\mu_{ G_{e}}\in \mathcal M_+(\mathcal T\times  G_{(j,k)}\times \Theta_j)$:
\red{\begin{align}
\mu_{G_{(j,k)}}(A\times B\times C)=\hspace*{-5mm}\int\limits_{\mathcal T\times D_j\times \Theta_j}\hspace*{-5mm}\mu_{j}(A\times B\times C\mid \tau_k,x_0,\theta)\,d\mu_{s_{j}},
\label{eq:mu_G}
\end{align}}%
for all $(j,k) \in {\mathcal E}$ are defined. The measures $\mu_{ G_{(j,k)}}$ are supported on the guards of mode $j$ and should be interpreted as the hitting times of the guard. \red{Note that $\mu_{G_{(i,j)}}$ has been defined with a slight abuse of notation. Guards were introduced in Defn.~\ref{def:system} as subsets of $\partial D$. In the uncertain description, if the guard description does not depend on $\theta$, then $\spt(\mu_{G_e})=G_{(j,k)}\times \Theta_j$; otherwise $\spt(\mu_{G_e})\subset G_{(j,k)}\times \Theta_j$.}
\begin{rem}
\label{rem:mu_G:restriction}
  Note the similarity between the definitions of $\mu_\mathcal J$ and $\mu_{G_\mathcal E}$ (Eqns.~(\ref{eq:mu_avg}) and (\ref{eq:mu_G})); $\mu_{G_{(j,k)}}$ is the restriction of $\mu_j$ to $G_{(j,k)}$ (denoted by $\mu_{j}|_{G_{(j,k)}}$).
\end{rem}
The {\em final measure} in each mode $j$ is then defined as:
\begin{align}
  \mu_{f_j}=\delta_T\otimes \mu_{T_j}+\sum_{k\in\{l\mid (j,l)\in \mathcal E\}}\mu_{ G_{(j,k)}}.
\label{eq:mu_T}
\end{align}
To compute $X_0$, we relate $\{\mu_{s_j}\}_{j\in \mathcal J}$ with $\{\mu_{f_j}\}_{j\in \mathcal J}$ using the dynamics of the system.
As a first step, define linear operators {$\mathcal L_{ f_j}\colon \mathcal C^1(\mathcal T\times D_j\times \Theta_j)\rightarrow \mathcal C(\mathcal T\times D_j\times \Theta_j)$} as:
\begin{align}
      \mathcal L_{f_j}v=\frac{\partial v}{\partial t}+\ip{\nabla_x v,\tilde f_j}
    \label{eq:Lv}
\end{align}
where $v\in \mathcal C^1(\mathcal T\times D_j\times \Theta_j;\R)$ is an arbitrary test function and $\nabla_x v$ computes the gradient of $v$ in the $D_j$ coordinates.
Suppose the system transitioned to mode $j$ at \mbox{$t=\tau_{k-1}$} with the state vector taking value upon reset $x(\tau_{k-1})$ and $\theta$.
The value of $v$, evaluated along the flow of the system and at $t=\tau_{k}$ is computed using the Fundamental Theorem of Calculus:
\begin{align}
\begin{aligned}
    v& \big(\tau_k,x(\tau_{k}\mid x(\tau_{k-1}),\theta_{k-1})\big)=v(\tau_{k-1},x(\tau_{k-1}),\theta_{k-1})\\
    &+\int_{\tau_{k-1}}^{\tau_{k}}\hspace*{-2mm}\mathcal L_{f}v(t,x(t\mid \tau_{k-1},x(\tau_{k-1}),\theta_{k-1}))\,dt.
\end{aligned}
\label{eq:FTC}
\end{align}
Using Eqn.~(\ref{eq:mu_lambda}), Eqn.~(\ref{eq:FTC}) can be re-written as:
\begin{align}
\begin{aligned}
    v\big(\tau_k,&x(\tau_{k}\mid \tau_{k-1},x(\tau_{k-1}),\theta_{k-1})\big)=v(\tau_{k-1},x(\tau_{k-1}),\theta_{k-1})\\
    &+\ip{\mu_j(\cdot\mid \tau_{k-1},x(\tau_{k-1},\theta_{k-1}),\mathcal L_{ f}v},
\end{aligned}
\end{align}
which can be simplified further by using Eqns.~(\ref{eq:mu_avg})--(\ref{eq:mu_T}):
\begin{align}
  \ip{\mu_{f_j},v}=\ip{\mu_{s_j},v}+\ip{\mu_{j},\mathcal L_{f}v}.
  \label{eq:liouville_1}
\end{align}
Alternatively, using the standard definition of adjoint operators\footnote{A linear operator $\mathcal L$ and its adjoint, $\mathcal L'$, satisfy the following relation:
$$    \ip{\mathcal L'\mu,v}=\ip{\mu,\mathcal Lv}.$$
}, Eqn.~(\ref{eq:liouville_1}) is re-written as:
\begin{align}
\ip{\mu_{f_j},v}=\ip{\mu_{s_j},v}+\ip{\mathcal L'_{f}\mu_{j},v}.
  \label{eq:liouville_2}
\end{align}
\emph{Eqn.~(\ref{eq:liouville_2}) defines a linear relation that initial and final measures evolving according to the hybrid dynamics must satisfy.} \red{Lemma~\ref{lemma:existence} formalizes this relation between trajectories of the system and Eqn.~(\ref{eq:liouville_2}).}
\par
With the average occupation and final measures defined, Lemma~\ref{lemma:initial_measure} completes the characterization of the Liouville eqn. by providing an explicit expression for the starting measure.
\begin{lem}
\label{lemma:initial_measure}
  The initial measure in each mode $j$, $\mu_{s_j}$, is expressible as the following
  \begin{align}
    \mu_{s_j}=\delta_{0}\otimes \mu_{0_j}\otimes \mu_{\theta_j}+\hspace*{-2mm}\sum_{i\in \{k\mid (k,j)\in \mathcal E\}}\hspace*{-3mm} R_{(i,j)}^*(\pi_{(t,x)}^*\mu_{ G_{(i,j)}})\otimes \mu_{\theta_j},
  \end{align}
  where $\pi_{(t,x)}^*$ denotes the pushforward constructed by lifting the $(t,x)$-projection operator, $\pi_{(t,x)}: \mathcal T\times D_j\times \Theta_j \to \mathcal T\times D_j$, to measures; and $R^*_{(i,j)}$ is the lifting of the reset map between modes $i$ and $j$.
\end{lem}
\begin{pf}
  During the execution of a hybrid system, any mode can be entered either at $t=0$ or due to a reset. By definition of quasi-uncertain systems (Defn.~\ref{def:system} and Alg.~\ref{alg:execution}), the value of the uncertain parameter in each mode is independent of the state's initial condition. Hence, the measure on initial states and uncertainty values is a product measure.
  \par
  Let $\delta_0\otimes \mu_{0_j}\times \mu_{\theta_j}\in \mathcal M_+(\{0\}\times D_j\times \Theta_j)$ be the measure on the set of initial conditions of system states at $t=0$; and $\sigma_{0_j}\in \mathcal M_+(\mathcal T\times D_j\times \Theta_j)$ be the measure on initial conditions because of trajectory resets. Then, the initial measure in the $(t,x)$-coordinate, using notations from Eqn.~(\ref{eq:mu_sj}), can be decomposed as:
\begin{align}
  \bar\mu_{0_j}=\delta_0\otimes\mu_{0_j}+\pi_{(t,x)}^*\sigma_{0_j}
\end{align}
where $\pi_{(t,x)}^*$ is the pushforward measure under the projection operator (refer to Chapter 11 in \cite{lee2003smooth} for an introduction to pushforwards).
\par
Trajectory resets occur if and only if trajectories reach any of the guards; hence $\sigma_{0_j}$ and $\mu_{G_e},\forall e\in \mathcal E$ must be related. To formalize this relationship notice that $\sigma_{0_j}$ can be decomposed into measures corresponding to the source of each reset:
\begin{align}
  \sigma_{0_j}=\sum_{i\in \{k\mid (k,j)\in \mathcal E\}} \sigma_{(i,j)}\otimes \mu_{\theta_j},
  \label{eq:sigma_def}
\end{align}
where $\sigma_{(i,j)}$ is the measure describing initial conditions that are reset into mode $j$ from guard $G_{(i,j)}$.
\par
Upon reaching guard $G_{(i,j)}$, the trajectories transition according to the reset map, $R_{(i,j)}$. Reset $R_{(i,j)}$ is independent of the value of the uncertainty (by Defn.~\ref{def:system}); and the uncertainty in the dynamics upon reset (if any), in mode $j$, is independent of the uncertainty (if any) in mode $i$. Thus, only the $(t,x)$-marginal of $\mu_{G_{(i,j)}}$ affects $\sigma_{(i,j)}$. Finally, since $\mu_{\theta_i}$ is a probability measure, the $(t,x)$-marginal is equal to $\pi_{(t,x)}^*\mu_{G_{(i,j)}}$.
Applying a change of variables formula we have for all $w\in \mathcal C(\mathcal T\times  D_j)$ \cite[Theorem 3.6.1]{Bogachev2007}:
\begin{align}
    \ip{\sigma_{(i,j)},w}=\ip{\pi_{(t,x)}^*\mu_{ G_{(i,j)}},w\circ R_{(i,j)}}.
    \label{eq:reset_measure}
\end{align}
Essentially, $\sigma_{(i,j)}$ is the pushforward measure of $\pi_{(t,x)}^*(\mu_{ G_{(i,j)}})$ under $R_{(i,j)}$, and we have the expression in the statement of the Lemma.
\end{pf}

\subsection{The primal}
\label{ssec:primal}

The problem of computing the uncertain backwards reachable set of ${\mathcal H}$ can be formulated as an infinite-dimensional linear program that supremizes the \emph{volume} of the set of initial condition:
  \begin{flalign}\nonumber
  & & \sup_{\Lambda} \hspace*{1cm} & \sum_{j \in {\mathcal J}}\ip{\mu_{0_j},\mathds{1} } && (P) \nonumber \\
  & & \text{st.} \hspace*{1cm} &\mu_{s_j}+\mathcal L_{f}'\mu_j=\,\mu_{f_j } && \forall j\in {\mathcal J} \label{eq:primal:liouville}\\
  & & & \mu_{0_j}+\hat\mu_{0_j}=\,\lambda_j && \forall j\in {\mathcal J} \\
  & & & \sum_{j \in {\mathcal J}}\ip{\mu_{T_j},\mathds{1}}=\,\sum_{j \in {\mathcal J}}\ip{\mu_{0_j},\mathds{1}} && \label{eq:mass_conservation}
  \end{flalign}
where $\lambda_j$ is the Lebesgue measure supported on $D_j$,\\ $\Lambda=\Big\{\big(\mu_{\mathcal J},\mu_{0_{\mathcal J}},\mu_{T_\mathcal J},\hat\mu_{0_{\mathcal J}},\mu_{G_{\mathcal E}}\big) \in \bigtimes\limits_{j \in {\mathcal J}} {\mathcal M}_+({\mathcal T} \times D_j \times \Theta_j) \bigtimes\limits_{j \in {\mathcal J}} {\mathcal M}_+( D_j )  \bigtimes\limits_{j \in {\mathcal J}} {\mathcal M}_+( X_{(T,j)}\times \Theta_j )$ $\bigtimes\limits_{j \in {\mathcal J}} {\mathcal M}_+( D_j ) \newline \bigtimes\limits_{e \in {\mathcal E}} {\mathcal M}_+( \mathcal T\times G_e\times \Theta_j)  \Big\}$, and $\mathds{1}$ denotes the function that takes value $1$ everywhere.
$\hat\mu_{0_j}\in \mathcal M(D_j)$ are slack variables introduced to ensure that the mass of the $\mu_{0_j}$ are identical to the volume (under the Lebesgue measure) of the uncertain backwards reachable set, as proven in Thm.~\ref{lemma:primal:main}.
Eqn.~(\ref{eq:mass_conservation}) ensures that all trajectories that emanate $\cup_{j\in \mathcal J}\,\spt(\mu_{0_j})$ reach $X_{T}$ at $t=T$.
\begin{thm}
\label{lemma:primal:main}
If $(\{\mu_{s_j}\}_{j \in {\mathcal J}},\{\mu_{j}\}_{j \in {\mathcal J}},\{\mu_{f_j}\}_{j \in {\mathcal J}})$ is a solution to $(P)$, then ${X}_{(0,j)} = \spt(\mu_{0_j})$ for each $j \in {\mathcal J}$.
In addition, the optimal value of $(P)$ is equal to the sum of \emph{volumes} of the uncertain backwards reachable set in each mode, i.e. $\sum_{j\in \mathcal J}\lambda_j(X_{0_j})$.
\end{thm}
\begin{pf}
We prove this Theorem by showing that the BRS contains the supports of $\mu_{0_j},\forall j\in \mathcal J$ and that the union of supports of $\mu_{0_j}$ contains the BRS.
\par
Suppose that $\bigcup_{j\in \mathcal J}\spt(\mu_{0_j})\backslash X_0 \ne \emptyset$; this implies that by Lemma~\ref{lemma:existence} (in the appendix), there exist trajectories that begin in $\bigcup_{j\in \mathcal J}(\spt(\mu_{0_j})\backslash X_{(0,j)})$ and terminate in $X_T$. This is a contradiction since the BRS is the largest set of initial conditions from which trajectories can reach $X_T$ at time $t=T$. Thus,
  \begin{align}
  &\bigcup_{j\in \mathcal J} \spt(\mu_{0_j})\subset \bigcup_{j\in \mathcal J} X_{(0,j)},\\
  &\sum_{j\in \mathcal J}\lambda_j(\spt(\mu_{0_j}))\le\sum_{j\in \mathcal J}\lambda_j( X_{(0,j)}).
  \label{eq:support_lemma:1}
  \end{align}
By definition of the BRS, all state trajectories that emanate from a subset of $X_0$ end in $X_T$. That is, for each $j\in \mathcal J$ and initial measure $\mu_{0_j}$, if $\spt(\mu_{0_j})\subset X_{(0,j)}$, there exist measures $\mu_{j}$ and $\mu_{f_j}$ that satisfy Eqn.~(\ref{eq:primal:liouville}).
Thus the following inequality is true:
    \begin{align}
    \sum_{j\in \mathcal J}\lambda_j(\spt(\mu_{0_j}))\ge \sum_{j\in \mathcal J}\lambda_j(X_{(0,j)})
      \label{eq:support_lemma:2}
    \end{align}
From Eqns.~(\ref{eq:support_lemma:1}) and (\ref{eq:support_lemma:2}), $\bigcup_{j\in \mathcal J}\,\spt(\mu_{0_j})$ is the BRS of the system.
That the optimal value of $(P)$ is the volume of the uncertain backward reachable set follows by noting that the slack variables ensure absolute continuity of each $\mu_{0_j}$ with respect to the Lebesgue measure and the observation that \mbox{$\lambda_{j}|_{X_{(0,j)}},\forall j\in \mathcal J$} is feasible in ($P$).
  \end{pf}
\subsection{The dual}
\label{ssec:dual}
The dual to $(P)$ for a quasi-uncertain hybrid system ${\mathcal H}$ can be written as:
\begin{flalign}
	 && \inf_F  \hspace*{1cm} &  \sum_{j\in {\mathcal J}}\ip{\lambda_j,w_j}                 && (D)\nonumber\\
	 && \text{st.} \hspace*{1cm} & w_j\ge \,0                                                 &&
      \begin{tabular}{@{}l@{}}$\forall x \in D_j$\\$\forall j\in \mathcal J$\end{tabular} \label{eq:dual:formulation:hatmu0}\\
	 &&&   v_j(T,x,\theta)+q\ge\, 0 ,\>                                  && \begin{tabular}{@{}l@{}} $\forall (x,\theta)\in \Phi_j$ \\ $\forall j\in \mathcal J$\end{tabular} \label{eq:dual:formulation:muT}\\
	 &&&  - \mathcal L_{f}v_j(t,x,\theta)\ge\,0 ,                       && \begin{tabular}{@{}l@{}} $\forall (t,x,\theta)\in \Omega_j$\\$\forall j\in \mathcal J$\end{tabular}  \label{eq:dual:formulation:muj}\\
	 &&&   w_j-\ip{\mu_{\theta_j},v_j(0,x,\theta)}-q\ge \,1,          && \begin{tabular}{@{}l@{}} $\forall x \in D_j$\\$\forall j\in \mathcal J$\end{tabular}  \label{eq:dual:formulation:mu0}\\
	 &&&   v_j -\ip{\mu_{\theta_k},v_k}\circ R_{(j,k)}\ge \,0,              && \begin{tabular}{@{}l@{}} $\forall (x,\theta)\in G_{(j,k)}\times  \Theta_j$\\ $\forall j\in \mathcal J$\\ $(j,k)\in \mathcal E$\end{tabular}\label{eq:dual:formulation:guard}
\end{flalign}
where $F = \Big\{ \big(v_{\mathcal J},w_{\mathcal J},q\big) \in \bigtimes\limits_{j \in {\mathcal J}}\mathcal C^1(\mathcal T\times D_j\times \Theta_j) \bigtimes\limits_{j \in {\mathcal J}} \mathcal C( D_j) \times \R \Big\}$, $\Phi_j = X_{(T,j)} \times \Theta_j $ and $\Omega_j = {\mathcal T} \times D_j \times \Theta_j$.
\begin{rem}
    The dual to $(P)$ as stated in $(D)$, is a result of a slight abuse of notations. By definition, $\forall (i,j) \in\mathcal E,\,R_{(i,j)}\colon D_i\rightarrow D_j$; however, for notational convenience, in $(D)$, $R_{(i,j)}$ is assumed to have been defined as $R_{(i,j)}\colon \mathcal T\times D_i\rightarrow \mathcal T\times D_j$ with the mapping being identity in the $t$-component and the standard $R_{(i,j)}$ in the $x$-components.
\end{rem}
The solution to $D$ can be used to determine the uncertain backwards reachable set:
\begin{lem}
	\label{lemma:dual_outerapprox}
If $\big(\{v\}_{j \in {\mathcal J}},\{w_j\}_{j \in {\mathcal J}},q\big)$ is a feasible point to $D$, then the super-level set:
\begin{align}
	\bigcup_{j\in \mathcal J}\,\{x \in D_j \mid w_j(x)\ge 1\}
\end{align}
is an outer approximation of the uncertain backwards reachable set of ${\mathcal H}$.
Furthermore there is a sequence of feasible solutions to $(D)$ such that for each $j \in {\mathcal J}$, the 1-super-level set of the feasible $w_j$ converges from above to the indicator function on $X_{(0,j)}$ in the $L^1$ norm and almost uniformly.
\end{lem}
\begin{pf}
To prove this lemma we project the uncertain backwards reachable set into each mode and show that it is part of the 1-level set of $w$.
Assume that the state trajectory terminates in $X_{(T,j_k)}$ for some $j_k$
The state trajectory must have arrived in mode $j_k$ through a finite sequence of mode-transitions (according to Assumption~\ref{assump:zeno}).
Let this sequences of mode-transitions be of length $k$.
Suppose the trajectory entered mode $j_k$ at time $\tau_k$, then from the Fundamental Theorem of Calculus (FTC) and the constraints in $(D)$, the following inequalities hold:
\begin{flalign}
	&&-q\le&\, v_{j_k}(T,x(T\mid x(\tau_k^+),\theta),\theta)&&\hspace*{-2mm} \text{\big(Eqn.~(\ref{eq:dual:formulation:muT})\big)}\\
    &&\le&\, v_{j_k}(\tau_k,x(\tau_k^+),\theta) && \hspace*{-2mm}\text{\big(FTC \& Eqn.~(\ref{eq:dual:formulation:muj})\big)}\label{eq:lemma:superlevel:FTC:end}
\end{flalign}
Now, integrating both sides of Eqn.~(\ref{eq:lemma:superlevel:FTC:end}) wrt. $\mu_{\theta_{j_k}}$, and noting that $\mu_{\theta_{j_k}}$ is a probability measure, we get
\begin{align}
	-q\le \ip{\mu_{\theta_{j_k}},v_{j_k}(\tau_k,x(\tau_k^+),\theta)}.
\end{align}
By iterative application of the constraints in $(D)$, it follows that:
\begin{flalign}
	&&-q\le&\, \ip{\mu_{\theta_{j_k}},v_{j_k}}\circ R_{(j_{k-1},j_{k})}(\tau_{k},x(\tau_{k}^-)) && \hspace*{-2mm}\text{(Reset)}\\
    &&  \le&\, v_{j_{k-1}}(\tau_{k},x(\tau_{k}^-\mid x(\tau_{k-1}^+),\theta),\theta) && \hspace*{-2mm} \text{\big(Eqn.~(\ref{eq:dual:formulation:guard})\big)}\\
    &&  \le &\, \ip{\mu_{\theta_{j_{k-1}}},v_{j_k}(\tau_k,x(\tau_{k-1}^+),\theta)} && \\\nonumber
    &&  &\,\vdots\\
    &&  \le &\,v_{j_0}(\tau_1,x(\tau_1^-\mid x_0,\theta),\theta)\\
    &&  \le &\,v_{j_0}(0,x_0,\theta)\\
     && \le &\,\ip{\mu_{\theta_{j_0}},v_{j_0}(0,x_0,\theta)}\\
     && \le &\, w_{j_0}(x_0)-q-1. && \hspace*{-2mm} \text{\big(Eqn.~(\ref{eq:dual:formulation:mu0})\big)}
\end{flalign}
The final inequality implies that the initial condition of every trajectory that ends in the terminal set belongs to the 1-superlevel set of $w_j$ for some $j \in {\mathcal J}$.
The remainder of the proof follows from a straightforward extension to \cite[Theorem 2]{henrion2014convex}.
\end{pf}
Finally, note that the value computed by either optimization problem is equal which follows from \cite[Theorem 3.10]{Anderson1987} and is similar to \cite[Theorem 2]{henrion2014convex}:
\begin{lem}
Formulations $(P)$ and $(D)$ are equivalent and have the same optimal value.
\end{lem}
\red{In the next section, we present a method to numerically solve the dual problem.}
\begin{rem}
There are two key aspects of the presentation that deserve re-iteration:
First, the uncertainties that influence the dynamics are drawn from the distribution each time a trajectory enters a new mode;
Second, the uncertain backwards reachable set corresponds to the set of initial conditions for \emph{all} trajectories that are able to reach the terminal set in spite of \emph{all} possible sequences of uncertainty that each have non-zero probability.
Notice that the uncertain backwards reachable set is the intersection of the backwards reachable set for every possible discrete uncertainty with non-zero probability.
\label{rem:prob_clarification}
\end{rem}

%% file: sections/sec_relax.tex
\section{Numerical Implementation}
\label{sec:implementation}
In this section, a sequence of Semidefinite Programs (SDP)s that approximate the solution to the infinite dimensional primal and dual defined in $\S$\ref{ssec:primal} and $\S$\ref{ssec:dual} are introduced.
This sequence of relaxations is constructed by characterizing each measure using a sequences of moments\footnote{The $n$th moment of a measure ($\mu$) is obtained by evaluating the following expression
  $$y_{\mu,n}=\ip{\mu,x^n}.$$}
and assuming the following:
\begin{assum}
The vector field in each mode and reset map between modes is a polynomial.
Moreover the domain, the value of uncertainties, the guard, and the target set in each mode is a semi-algebraic set.
  \label{assump:poly}
\end{assum}
Recall that polynomials are dense in the set of continuous functions by the Stone-Weierstrass Theorem \cite{folland2013real}; so this assumption is made without much loss of generality.
\par
Under this assumption, given any finite $d$-degree truncation of the moment sequence of all measures in the primal $(P)$, a primal relaxation, $(P_d)$, can be formulated over the moments of measures to construct an SDP.
The dual to $(P_d)$, $(D_d)$, can be expressed as a sums-of-squares (SOS) program by considering $d$-degree polynomials in place of the continuous variables in $D$.
\par
To formalize this dual program, first note that a polynomial $p \in \R[x]$ is SOS or $p \in \text{SOS}$ if it can be written as $p(x) = \sum_{i=1}^m q_i^2(x)$ for a set of polynomials $\{q_i\}_{i=1}^m \subset \R[x]$.
Note efficient tools exist to check whether a finite dimensional polynomial is SOS using SDPs~\cite{parrilo2000structured}.
Next, suppose we are given a semi-algebraic set $A = \{x \in \R^n \mid h_{i}(x) \geq 0, h_i \in \R[x], \forall i \in \N_m \}$.
We denote the $d$-degree {\em quadratic module} of $A$ as:
\begin{align}
  \begin{split}
  Q_d(A)=\bigg\{q\in \R_d[x]\,\bigg|\, \exists \{s_k\}_{k \in N_m \cup \{0\}} \subset \text{SOS s.t. } q=s_0+\sum_{k\in \N_{m}}h_{k}s_k \bigg\}
  \end{split}
\end{align}

\noindent The $d$-degree relaxation of the dual, $D_d$, can be written as:
  \begin{flalign}
    & & \inf_{\Xi_d} \hspace*{0.1cm} & \sum_{j\in \mathcal J}\int_{D_j}w_j(x)\,d\lambda_j(x) && \hspace*{-0.3cm} (D_d) \\
    & & \text{st.} \hspace*{0.1cm} & w_j^d\in Q_d(X_{(T,j)}) && \hspace*{-0.3cm} \forall j\in \mathcal J \\
    & & & \hspace*{-3mm}v_j^d(T,\cdot)+q\in Q_d(D_{j} \times \Theta_j) && \hspace*{-0.3cm} \forall j\in \mathcal J \\
    & & & \hspace*{-3mm}-\mathcal L_{f_j}v_j^d\in Q_d(\mathcal T\times D_{j} \times \Theta_j) && \hspace*{-0.3cm} \forall j\in \mathcal J \\
    & & & \hspace*{-3mm}w_j^d-\ip{\mu_{\theta_j},v_j^d(0,\cdot)}-q-1\in Q_d(D_j) && \hspace*{-0.3cm} \forall j\in \mathcal J \\
    & & & \hspace*{-3mm}v_j^d-\ip{\mu_{\theta_k},v_k^d}\circ R_{(j,k)}\in Q_d(\mathcal T \times D_j \times \Theta_j) && \hspace*{-0.3cm} \forall (j,k)\in {\mathcal E}
  \end{flalign}
where $\Xi_d=\Big\{ \big(v^d_{\mathcal J},w^d_{\mathcal J},q \big) \in \bigtimes\limits_{j\in {\mathcal J}}\R_d[t,x,\theta]$ $\bigtimes\limits_{j \in {\mathcal J}} \R_d[x]\times\R \Big\}$.
A primal can similarly be constructed, but the solution to the dual can be used directly generate a sequence of outer approximations to the uncertain backwards reachable set:
\begin{thm}
	For each $d \in \N$ and $j \in {\mathcal J}$, let $w_{j_d}$ denote the $j$-slice of the $w$-component of the solution to $D_d$. Then ${ X}_{(0,j_d)} = \{x \in D_j \mid w_{j_d}(x) \geq 1 \}$ is an outer approximation to ${ X}_{(0,j)}$ and $\lim_{d\to\infty}\lambda_{n_j}({ X}_{(0,j_d)} \backslash { X}_{(0,j)}) = 0$.
\end{thm}
\begin{pf}
  The proof to this lemma is an extension of Theorems 5--7 in \cite{shia2014convex} given Lemma~\ref{lemma:dual_outerapprox}.
\end{pf}

%% file: sections/sec_extensions.tex
\section{Extensions}
\label{sec:ext}
In this section, we present extensions to the problem formulation presented in $\S$\ref{sec:prob}. All results presented in $\S$\ref{sec:prob} and \ref{sec:implementation} hold for the formulations contained herein; proofs are largely identical and are omitted.
\subsection{Free terminal time}
\label{ssec:free_time}
The formulation in $\S$\ref{sec:prob} aimed at identifying the BRS of quasi-uncertain hybrid systems given a fixed terminal time $T$ and a terminal set $X_T$. Suppose it is of interest to estimate the BRS which includes all initial conditions from which trajectories reach $X_T$ at some time $t\le T$. This problem of computing the uncertain backwards reachable set of ${\mathcal H}$ can be formulated as the following infinite-dimensional linear program that supremizes the \emph{volume} of the set of initial conditions:
  \begin{flalign}\nonumber
  & & \sup_{\Lambda} \hspace*{1cm} & \sum_{j \in {\mathcal J}}\ip{\mu_{0_j},\mathds 1} && (P^T) \nonumber \\
  & & \text{st.} \hspace*{1cm} &\mu_{s_j}+\mathcal L_{f}'\mu_j=\,\mu_{f_j } && \forall j\in {\mathcal J} \\
  & & & \mu_{0_j}+\hat\mu_{0,j}=\,\lambda_j && \forall j\in {\mathcal J} \label{eq:free_terminal_w}\\
  & & & \sum_{j \in {\mathcal J}}\ip{\mu_{T_j},\mathds 1}=\,\sum_{j \in {\mathcal J}}\ip{\mu_{0_j},\mathds 1} &&
  \end{flalign}
  where $\lambda_j$ is the Lebesgue measure supported on $D_j$,\\ $\Lambda=\Big\{\big(\mu_{\mathcal J},\mu_{0_{\mathcal J}},\mu_{T_\mathcal J},\hat\mu_{0_{\mathcal J}},\mu_{G_{\mathcal E}}\big) \in \bigtimes\limits_{j \in {\mathcal J}} {\mathcal M}_+({\mathcal T} \times D_j \times \Theta_j) \bigtimes\limits_{j \in {\mathcal J}} {\mathcal M}_+( D_j )  \bigtimes\limits_{j \in {\mathcal J}} {\mathcal M}_+( X_{(T,j)}\times \Theta_j )$ $\bigtimes\limits_{j \in {\mathcal J}} {\mathcal M}_+( D_j ) \newline \bigtimes\limits_{e:=(j,k) \in {\mathcal E}} {\mathcal M}_+( \mathcal T\times G_e\times \Theta_j)  \Big\}$.
  \par
  \red{This formulation differs from $(P)$ only in that $\mu_{T_j}$, is now supported on $\mathcal T\times X_{(T,j)}\times \Theta_j$ as opposed to $X_{(T,j)}\times \Theta_j$. This change in support translates into admitting initial conditions from which solution trajectories reach $X_T$ at some time before $t=T^+$ for all sequences of uncertainties, as a part of the BRS. Note that this does \emph{not} mean that solution trajectories remain in $X_{(T,j)}\times \Theta_j$ for all time, upon first entry.}
\subsection{\red{Inner approximations}}
As as review, according to the execution of quasi-uncertain hybrid systems (Alg.~\ref{alg:execution}), the BRS of a set $X_T$ is the set of initial conditions that reach $X_T$ for all possible sequences of uncertainties (refer to $\S$\ref{sec:preliminaries} for definition). Given a terminal set $X_T$, the problem formulation presented in $\S$\ref{sec:prob} and the subsequent relaxations in $\S$\ref{sec:implementation} provide outer approximations of the BRS. In some situations, inner approximations of the BRS are more informative since they guarantee that points in the interior of the set obtained through relaxations will satisfy the problem objectives. In this section, we present a problem formulation whose relaxations provide convergent inner approximations of the BRS by adapting the technique presented in \cite{Korda2013}.
\par
For {quasi-uncertain hybrid} systems with the BRS as defined in $\S$\ref{sec:preliminaries}, constructing a sequence of convergent \emph{inner approximations} of the BRS of a set $X_T$ is analogous to constructing a sequence of convergent \emph{outer approximations} of the set of initial conditions for which there exists a set (of non-zero measure) of sequences of uncertainties, $(\theta)_n$, for which the resulting trajectories fail to arrive at the terminal set at $t=T$.
\par
The infinite-dimensional problem of interest in this case, employing the nomenclature adopted in $\S$\ref{sec:prob}, is the following:
  \begin{flalign}\nonumber
  & & \sup_{\Lambda} \hspace*{1cm} & \sum_{j \in {\mathcal J}}\ip{\mu_{0_j},\mathds 1} && (P^p) \nonumber \\
  & & \text{st.} \hspace*{1cm} &\mu_{s_j}+\mathcal L_{f}'\mu_j=\,\mu_{f_j } && \forall j\in {\mathcal J} \\
  & & & \mu_{0_j}+\hat\mu_{0,j}=\,\lambda_j && \forall j\in {\mathcal J} \label{eq:aggregate:w}\\
  & & & \sum_{j \in {\mathcal J}}\ip{\mu_{T_j},\mathds 1}=\,\sum_{j \in {\mathcal J}}\ip{\mu_{0_j},\mathds 1} &&
  \end{flalign}
  where $\lambda_j$ is the Lebesgue measure supported on $D_j$,\\ $\Lambda=\Big\{\big(\mu_{\mathcal J},\mu_{0_{\mathcal J}},\mu_{T_\mathcal J},\hat\mu_{0_{\mathcal J}},\mu_{G_{\mathcal E}}\big) \in \bigtimes\limits_{j \in {\mathcal J}} {\mathcal M}_+({\mathcal T} \times D_j \times \Theta_j) \bigtimes\limits_{j \in {\mathcal J}} {\mathcal M}_+( D_j )  \bigtimes\limits_{j \in {\mathcal J}} {\mathcal M}_+\big( (X_{(T,j)}^c\bigcap D_j)\big)$ $\bigtimes\limits_{j \in {\mathcal J}} {\mathcal M}_+( D_j ) \newline \bigtimes\limits_{e:=(j,k) \in {\mathcal E}} {\mathcal M}_+( \mathcal T\times G_e\times \Theta_j)  \Big\}$.
\par
Problem $(P^p)$ differs from $(P)$ on two main counts -- (1) how uncertainty in incorporated; (2) definition of the final measure in each mode. In problem $(P)$, uncertainty was augmented as a state with a given initial distribution; in $(P^p)$, the uncertainty is considered as a bounded input to the system. This difference harks back to the objective of the inner approximation problem -- find all initial conditions that fail to reach $X_T$ for \emph{some} sequence of uncertainty. By assuming that $\theta$ is a bounded control whose value can be arbitrarily chosen, we are able to search over all possible uncertainty sequences (which includes, specifically, the class of uncertainties that is constant in each mode and changes value upon trajectory reset).
\par
For solution trajectories to fail to arrive ar $X_T$ at $t=T$, one of the following must be true: solution trajectories either leave the space at some time\footnote{For a hybrid system, given that guards are subsets of the boundary (Remark~\ref{rem:guard_boundary}), leaving the space is equivalent to arriving at $\bigcup_j(\partial D_j - \bigcup_{k\in\{l\mid (j,l)\in \mathcal E\}}G_{(j,k)})$ at $t\in [0,T]$.} $t\in [0, T]$ or they arrive at $D\backslash X_T$ at $t=T$ for some sequence of uncertainty. That is, employing the same nomenclature as in $\S$\ref{sec:prob}, the \emph{final measure} in each mode is give by
\begin{align}
  \mu_{f_j}=\delta_T\otimes \mu_{T_j}+\mu_{\partial_j}+\sum_{k\in\{l\mid (j,l)\in \mathcal E\}}\mu_{ G_{(j,k)}}
\end{align}
where $\mu_{T_j}\in \mathcal M_+\big((X_{(T,j)}^c\bigcap D_j)\times \Theta_j\big),\,\mu_{\partial_j}\in \mathcal M_+\big(\mathcal T\times (\partial D_j-\bigcup_{k\in\{l\mid (j,l)\in \mathcal E\}}G_{(j,k)})\big)$ and $\mu_{G_{(j,k)}}\in \mathcal M_+(\mathcal T\times G_{(j,k)}\times \Theta_j),\,\forall (j,k)\in \mathcal E$.
\par
In the above definition, $\mu_{\partial _j}$ is a measure that \emph{traps} trajectories that leave the domain, $D_j$, and not pass through any of the guards. Note that $\mu_{T_j}$ in $(P^p)$, unlike in $(P)$, is supported on the complement of $X_{(T,j)}$; that is, mirroring out desire to characterize the set of initial conditions from which solution trajectories fail to reach $X_{(T,j)}$.
\par
It should be remarked that to be able to derive guarantees/results for $(P^p)$ as we have done for $(P)$, we need the following technical assumption.
\begin{assum}
$X_T$ is an open subset of $D$; hence $X_T^c$ is closed; and  $X_T^c\bigcap (\bigcup_{e\in \mathcal E}G_{e})=\emptyset$.
\end{assum}

%% file: sections/sec_examples.tex
\section{Examples}
\label{sec:examples}
In this section, the proposed method is applied to six examples of varying complexities. The examples are chosen to highlight different variants of the problem description:
\begin{enumerate}
  \item The compass-gait without uncertainty shows the applicability to deterministic systems
  \item The 1D example in Ex.~\ref{example:1D} and another example adopted from \cite{Chesi2004} are used to as benchmark examples for computing the uncertain BRS for linear and nonlinear systems
  \item The rimless wheel is an example in which guards are uncertain
  \item The bean-bag toss serves as an example for the free terminal-time formulation in $\S$\ref{ssec:free_time} and frictional cone uncertainties
  \item A variant of logistic resource growth dynamics is used to demonstrate inner approximations
\end{enumerate}
The relaxed dual problems are constructed using the SPOTLESS toolbox \cite{spotless} and solved with MOSEK on a computer equipped with a Intel Xeon W3540 processor and 12GB of RAM. Note that, considering the definition of BRS in $\S$\ref{sec:preliminaries}, for all probability distributions with identical support, the problem formulation will generate the same uncertain backwards reachable set. As a result it is assumed that the disturbance, $\theta$, is uniformly distributed (denoted as $\theta\sim\mathcal U([a,b])$ for $\theta$ uniformly distributed in the interval $[a,b]$). Additionally for numerical stability, the domains of modes of the hybrid system are scaled to a box of the appropriate dimension.
\subsection{Compass Gait Walker without Uncertainty}
\begin{figure}[!t]
\centering
  {\includegraphics[width=0.5\columnwidth,trim=0cm 0cm 0cm 0cm, clip=true]{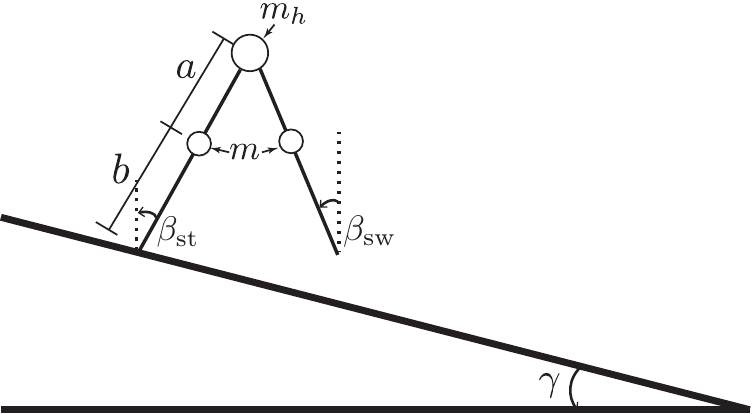}}
  \caption{Schematic of the compass gait walker.}
    \label{fig:cg_model}
\end{figure}

The compass gait (CG) walker is a simple model of legged locomotion consisting of two legs: one leg fixed to the ground called the stance leg, and one leg that swings called the swing leg and shown in Fig.~\ref{fig:cg_model}.
The CG is a one mode hybrid system in which the guard is reached when the swing leg makes contact with the inclined plane, upon which the swing and stance leg switch.
With different slopes and parameters, the CG has been found to reach a limit cycle consisting of 1, 2+ steps \cite{Goswami1998}.
In this example, we consider a passive CG walker with no actuation.
Let $\beta=[\beta_{\text{sw}}, \beta_{\text{st}}]$ and $l=a+b$, the dynamics of the passive CG are given by:
\begin{equation}
M(\beta,\dot{\beta})\ddot{\beta} + C(\beta,\dot{\beta})\dot{\beta} + N(\beta) = 0\end{equation}
where
\begin{equation}
M(\beta,\dot{\beta}) = \left[ \begin{array}{cc}
mb^2 & -mlb\cos(\beta_{\text{st}}-\beta_{\text{sw}}) \\
-mlb\cos(\beta_{\text{st}}-\beta_{\text{sw}}) & (m_h + m)l^2  + ma^2
\end{array} \right]
\label{eq:cg_mass}
\end{equation}

\begin{equation}
C(\beta,\dot{\beta}) = \left[ \begin{array}{cc}
0 & mlb\sin(\beta_{\text{st}}-\beta_{\text{sw}})\dot{\beta}_{\text{st}} \\
mlb\sin(\beta_{\text{st}}-\beta_{\text{sw}})\dot{\beta}_{\text{sw}} & 0
\end{array} \right]
\label{eq:cg_coriolis}
\end{equation}

\begin{equation}
N(\beta) = \left[ \begin{array}{c}
mbg \sin(\beta_{\text{sw}}) \\
-(m_h l + ma + ml)g \sin(\beta_{\text{st}})
\end{array} \right]
\label{eq:cg_potential}
\end{equation}

The guard is defined as when the swing leg hits the inclined slope and mathematically defined as:
\begin{align}
G_{(1,1)}=\{ (\beta, \dot{\beta})\mid \beta_{\text{sw}} + \beta_{\text{st}}+ 2\gamma = 0\}.
\end{align}

The reset map is given by:
\begin{align}
R_{(1,1)}(\beta^-,\dot \beta^-) = \begin{bmatrix}
  \beta_{\text{st}} &
  \beta_{\text{sw}} &
  \left(Q_\alpha^+\right)^{-1}Q_\alpha^- \dot{\beta}^-
  \end{bmatrix}'
\end{align}
where
\begin{equation}
Q_\alpha^- = \left[ \begin{array}{cc}
-mab & -mab + (m_h l^2 + 2mal)\cos(2 \alpha) \\
0 & -mab
\end{array} \right]
\end{equation}
\small
\begin{equation}
Q_\alpha^+ = \left[ \begin{array}{cc}
mb(b-l\cos(2\alpha)) & ml(l - b \cos(2\alpha)) + ma^2 + m_h l^2 \\
mb^2 & -mbl\cos(2\alpha)
\end{array} \right]
\end{equation}
\normalsize

and $\alpha = \frac{\beta_{\text{sw}}-\beta_{\text{st}}}{2}$.
The reset dynamics are derived using conservation of momentum resulting in a loss of kinetic energy.
The loss of kinetic energy is recovered via change in potential energy as the CG walks down the slope.


Prior to \cite{Manchester2011} which presents an inner-approximation to the BRS, the BRS was limited to exhaustive simulation.
However, \cite{Manchester2011} is limited to a small region and misses much of the BRS.
For computation, we consider the 5th order Taylor approximation of the dynamics about the origin and a linearized reset map about the point where the limit cycle encounters the guard.

\begin{figure}[!t]
\centering
  {\includegraphics[width=0.7\columnwidth,trim=0cm 0cm 0cm 0cm, clip=true]{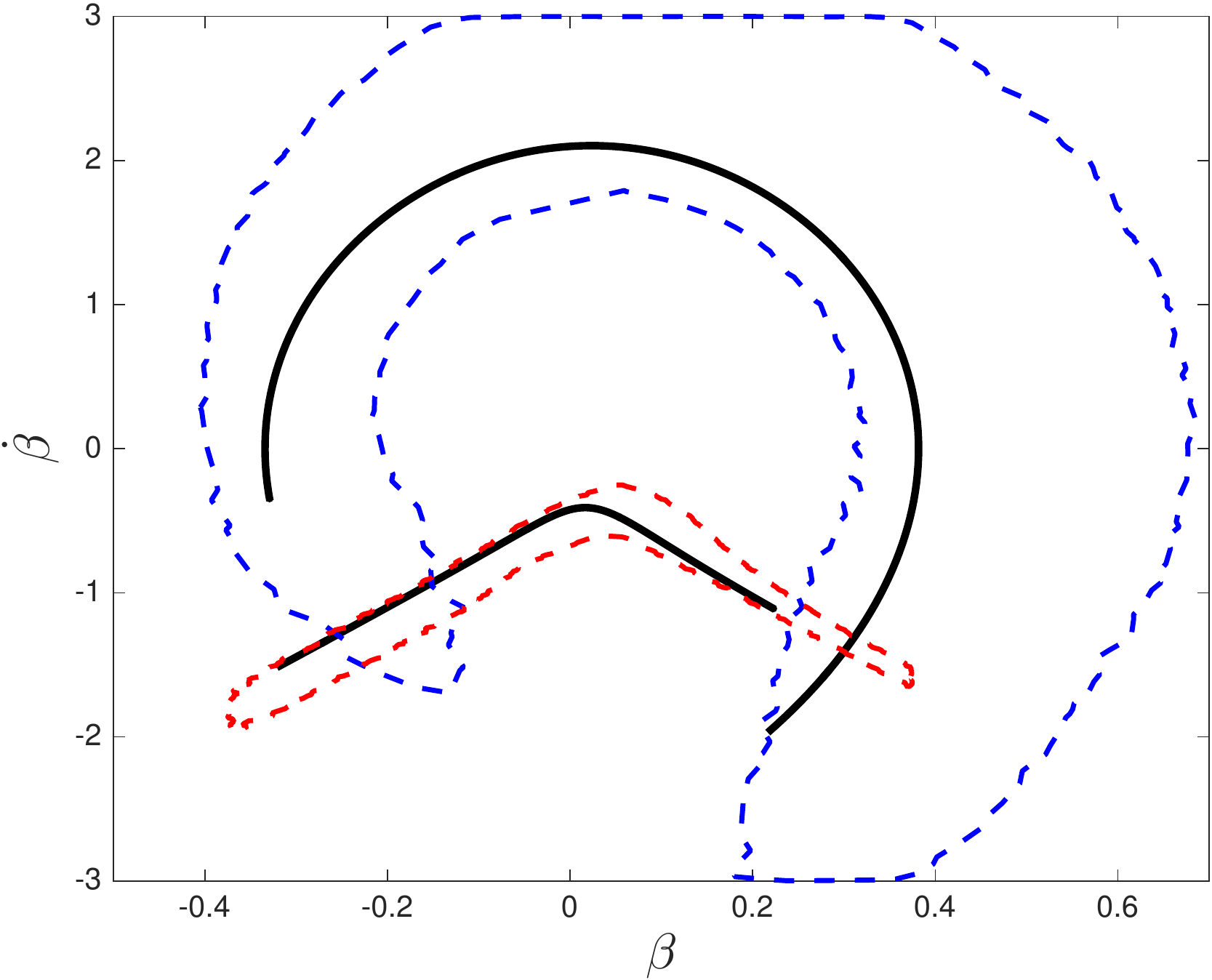}}
  \caption{The finite time region of attraction, shaded in gray, for the compass gait's limit cycle. The swing and stance leg limit cycles are projected down to the ($\beta,\dot\beta$) domain and denoted by the blue and red line, respectively. The time horizon is 1s.}
    \label{fig:cg_brs}
\end{figure}

Figure~\ref{fig:cg_brs} presents the polynomial degree 10 approximation to the backwards reachable set for the compass gait (with $\gamma = 0.05$) which is tasked with reaching within $0.1$ of the limit cycle (in black) in $T=1$ second with $m_h = 10$kg, $m=5$kg, $a=b=1$.
Through simulation of 10,000 points in the BRS, we find that 70\% of the BRS reaches within 0.1 of the limit cycle.

\subsection{1-D Quasi-Uncertain Linear System}
\begin{figure}[!t]
\centering
  {\includegraphics[width=0.7\columnwidth,trim =1.5in 4.5in 1.5in 3.5in, clip=true]{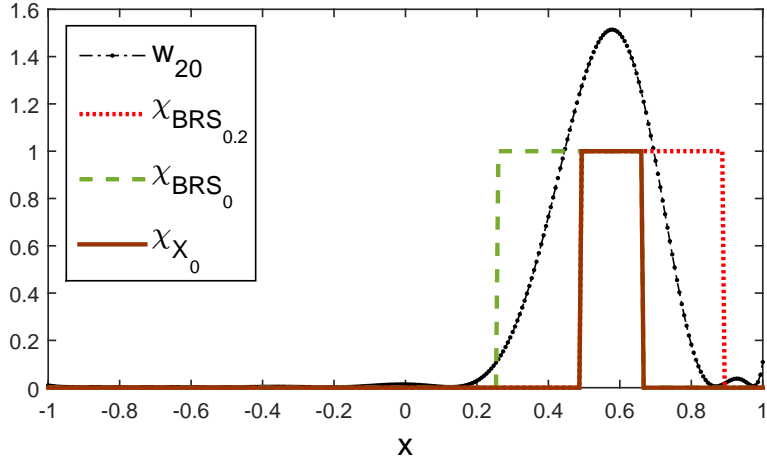}}
  \caption{Comparison between indicator functions on analytically computed backwards reachable sets of the 1D linear system in Eqn.~(\ref{eq:ex:1d}) and $w_{20}$ obtained as a part of the solution to $D_{20}$. The 1-level superset of $w_{20}$ encompasses the intersection of the backward reachable sets when the uncertainty takes extreme constant values.}
    \label{fig:1D:linear}
\end{figure}
Recall the 1-D linear dynamical system from Ex.~\ref{example:1D} whose dynamics are:
  \begin{align}
	\dot x_1 = -0.7x_1+0.2\theta-0.1,
\label{eq:ex:1d}
\end{align}
where $\theta\in \mathcal U([0.2,1])$.
Setting $T = 1$, the target set is chosen as $X_T=[0.2,0.4]$.
If $\theta$ was a fixed constant then the BRS for the system evolving with this known constant is analytically computed to be
\begin{align}
	\begin{aligned}
  		BRS_\theta=&\left[\left(0.2-\frac{2\theta-1}{7}\right)e^{0.7}+\frac{2\theta-1}{7}, \right. \\
					&\hspace*{0.2cm}\left.\left(0.4-\frac{2\theta-1}{7}\right)e^{0.7}+\frac{2\theta-1}{7}\right].
	\end{aligned}
\end{align}
\normalsize
Note that the expression for the $BRS_\theta$ is linear in $\theta$ and that the width of $BRS_\theta$ is constant for all values of $\theta$.
As the value of $\theta$ changes, $BRS_\theta$ slides along $\R$; thus, the intersection of $BRS_{0.2}$ and $BRS_1$ is the uncertain backwards reachable set of the system in Eqn.~(\ref{eq:ex:1d}) system.
\par
Figure~\ref{fig:1D:linear} plots the degree 20 approximation of the indicator function of $X_0$, $w_{20}$, that solves $D_{20}$ and the analytically computed indicator functions supported on $BRS_{0.2}$ and $BRS_{1}$, and $\chi_{X_0}$. Observe that the 1-level superset of $w_{20}$ contains $X_{0}$ and hence is an outer approximation of $X_0$.
\subsection{Rimless Wheel on Uneven Terrain}
\begin{figure}[!t]
\centering
  \includegraphics[trim=1.5in 3.3in 1.5in 3.5in, clip=true,width=0.7\columnwidth]{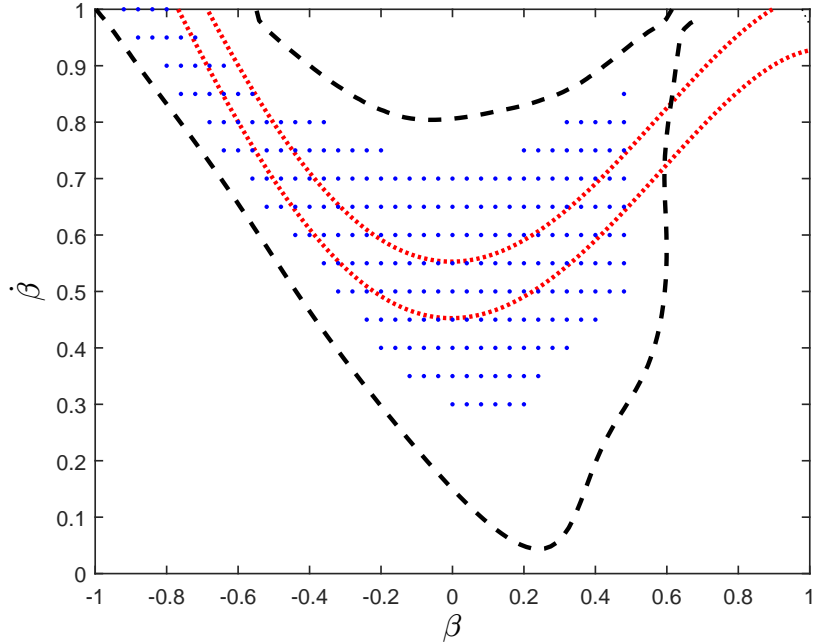}
  \caption{Comparing outer approximations with Monte Carlo simulations. The dotted curve depicts the boundary of the outer approximation of the backward reachable set of a PRW tasked with reaching the red band by $t=4$ seconds. The blue dots represent a sampling of the {\em true} backwards reachable set obtained by performing Monte Carlo simulation; for each initial condition, 100 trials were performed and only those initial conditions from which no trajectory violated constraints and reached $X_T$ at $t=4$ were included.}
  \label{fig:rw_brs}
\end{figure}
The rimless wheel, introduced in Ex.~\ref{example:rw}, is a one mode hybrid system in which the guard is reached when the swinging spoke makes contact with the inclined plane. For a rimless wheel rolling along an inclined plane with no terrain height variation (apart from the deterministic incline), an analytically computable stable limit cycle exits \cite{Coleman1998}; however, for the case considered in this example---with the inclined plane having variations in terrain height---the definition of a limit cycle is less clear. In this example, we define the terminal set as an $\epsilon$ band around the stable limit-cycle of the disturbance-free system.
\par
Figure~\ref{fig:rw_brs} presents the polynomial degree 14 approximation to the uncertain backwards reachable set (black dashed) for the rimless wheel (with $\alpha=0.4$) which is tasked with arriving within the red band in $T=4$ seconds, as it is rolling down an inclined plane with slope $\gamma=0.2$. The uncertain parameter, $\theta$, which affects the terrain height as described earlier in Ex.~\ref{example:rw}, is drawn from a uniform distribution, $\theta\sim\mathcal U([-0.1,0.1])$. The maximum terrain variation is about 25\% of the length of each spoke. With this setup, by the terminal time, somewhere between four and six spokes will have made contact with the wedge.
\par
The outer approximation of the uncertain BRS is validated by performing Monte Carlo simulations; the unit box is discretized into 51 points in either direction and 100 independent trajectories are simulated (using MATLAB's {\em ode45} function) from each initial condition.
The blue dots depict the initial conditions that arrived within the terminal set at the desired time without violating any constraint. Note that the set of points that succeeded in the Monte Carlo simulation is contained entirely in the uncertain BRS computed using our formulation. In fact, as a result of our method, we know that for points outside of the black region there exist a sequence of terrain heights that produces a trajectory that does not arrive at the target set at the designated time.
\par
Monte Carlo simulations for this example were performed on a computer with two Intel Xeon E5-2660 processors and 128 GB of RAM. The program was written and executed in parallel in MATLAB 8.4.0 using the parallel processing toolbox and took 21019 seconds to compute the result. In comparison, the proposed method, when solved on an arguably lesser computational resource, took 4487 seconds to solve $D_{14}$.
\subsection{An example from literature}
\begin{figure}
\centering
\begin{tikzpicture}
\node at (.2,1.75) {\small $T=5$ s};
\node at (4.25,1.75) {\small $T=10$ s};
\node at (0,0) {{\includegraphics[width=0.4\columnwidth, trim=3.65in 2.6in 1.5in 5.8in ,clip=true]{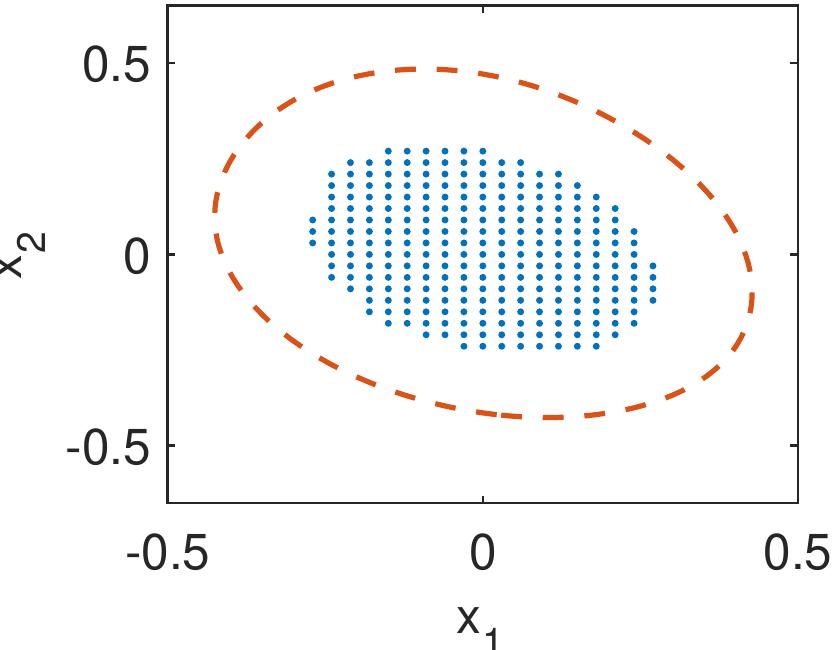}}};
\node at (7.5,0) {{\includegraphics[width=0.4\columnwidth, trim=1.2in 5.4in 4.0in 3.0in ,clip=true]{figures/examples/topcu_2010_ex1/BRS_T_5n10_p05_d_16_MC}}};
\end{tikzpicture}
\caption{Estimated degree-16 BRSs with terminal set $B_0(0.05)$ (dashed) and corresponding result of MC simulation (dots); (a) terminal time $T=5$ s, (b) terminal time $T=10$ s.}
\label{fig:topcu_chesi}
\end{figure}

We next consider the following benchmark example from \cite{Topcu2010,Chesi2004}
\begin{align}
\begin{split}
  \dot x =\begin{bmatrix}
    -x_1\\3x_1-2x_2
  \end{bmatrix}+\begin{bmatrix}
    -6x_2+x_2^2+x_1^3\\-10x_1+6x_2+x_1x_2
  \end{bmatrix}\delta+
  \begin{bmatrix}
    4x_2-x_2^2\\12x_2-4x_2
  \end{bmatrix}\delta^2
  \end{split}
\end{align}
with $\delta\sim \mathcal U([0,1])$.
We solve the outer approximation dual formulation with the terminal set described by a ball of radius 0.05 centered at the origin, for different values of terminal time, $T$. Figure~\ref{fig:topcu_chesi} presents a comparison between the obtained 16-degree outer approximations and the \emph{true} BRS as estimated by sampling the space for when $T=5$ s and $T=10$ s; in both cases, the 1-superlevel set is an outer approximation of the BRS.
\subsection{Bean-bag toss}
In this example, we consider a simplified version of the party game involving tossing a bean-bag onto an inclined plane which has a designated target-zone. If the motion of the bean-bag (of mass $m$) is considered as it leaves the hand of the player, its equation of motion can be modeled as a hybrid system with three states: (mode 1) flight, (mode 2) sliding, and (mode 3) rest.
The dynamic equations in each mode is given by the following
\begin{flalign}
  &\text{\bf Mode 1}&&\>\begin{bmatrix}
    \dot x_1& \dot x_2&
    \dot y_1& \dot y_2
  \end{bmatrix}^T =\, \begin{bmatrix}
    x_2&0&y_2&-g
  \end{bmatrix}^T&\\
  &\text{\bf Mode 2} &&\>\begin{bmatrix}
    \dot d_1\\\dot d_2
  \end{bmatrix} = \,\begin{bmatrix}
    d_2\\ mg\sin(\gamma) - sgn(d_2)\,\mu_kmg\cos(\gamma)
  \end{bmatrix}\\
  &\text{\bf Mode 3} &&\> \begin{bmatrix}
\dot d_1\end{bmatrix}= \,\begin{bmatrix}0\end{bmatrix}.
\end{flalign}
In the above description, in mode 1, the states correspond to position and velocities in cartesian coordinates; when the bean-bag makes contact with the wedge (of slope $\gamma$), it is assumed that the impact is inelastic and the system transitions to mode 2. Upon transitioning to mode 2, a coordinate transformation is applied and the states of the dynamic model correspond to the position of the bean-bag relative to the end of the wedge ($d$), and the sliding velocity along the wedge ($\dot d$). If the bean-bag slides along the wedge, it may reach the target-zone located at $d=0.55$ and shaped as a ball of radius $0.05$ m; whereupon the bean-bag falls and comes to rest (dynamics transitions to mode 3). It should be noted that as the bean-bag slides, its motion is impeded by friction with an uncertain kinetic frictional coefficient $\mu_k\in [0.37,0.97]$.
\par
The various guards and reset maps of this hybrid system are mathematically represented as
\begin{align}
  G_{(1,2)} = &\, \{(x_1,x_2,y_1,y_2)\in \R^4\mid y_1 = x_1\tan(\gamma)\},\\
  G_{(2,3)} = &\, \{(d_1,d_2)\in \R^2\mid d_2 =0\},\\
  R_{(1,2)} \mapsto &\, [y_1\csc(\gamma);x_2\cos(\gamma) + y_2\sin(\gamma)],\\
  R_{(2,3)} \mapsto &\, d_1.
\end{align}
Suppose the mass of the bean-bag is $m=0.5$ Kg, $\gamma =\pi/6$, and the length of the inclined face of the wedge is 1 m, and additionally that the velocity of the bean-bag as it slides $\dot d\in [-1,1]$. We are interested in knowing the positions, $(x_1,y_1)$, and velocities, $(x_2,y_2)$, at which one should toss the bag to ensure that it falls into the hole on the wedge.
\begin{figure}[!t]
\centering
  {\includegraphics[width=0.7\columnwidth,trim = 1.4in 3.3in 1.5in 3.5in, clip=true]{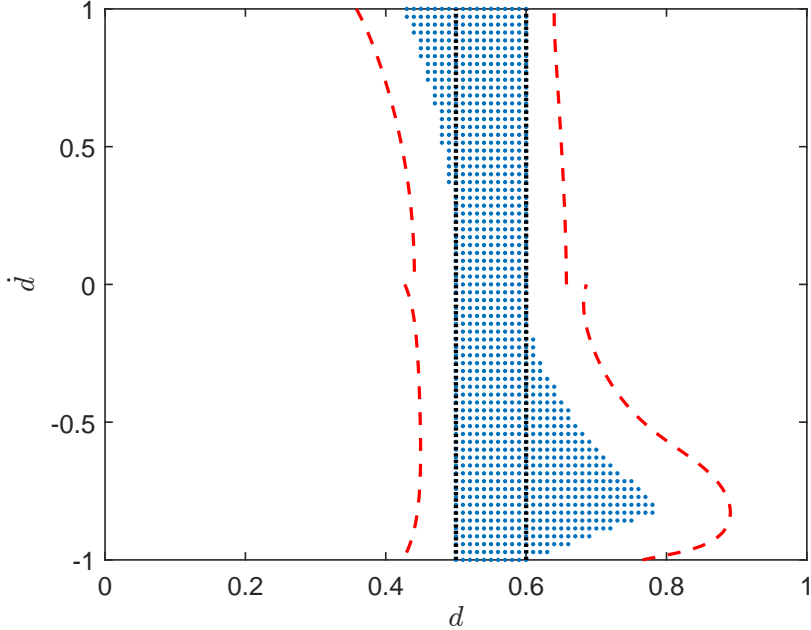}}
  \caption{Projection of the BRS of a bean-bag of mass 0.5 Kg thrown onto a 1 m long $\pi/6$-wedge with hole located at $X_T=[0.5,0.6]$. The kinetic frictional coefficient of the wedge's slope, $\mu_f\in[0.37,0.97]$, is not known with certainty. The blue dots are the landing points that hit the target; the interior of the red region in the outer approximation of the BRS.}
  \label{fig:bean_bag_toss}
\end{figure}
In Fig.~\ref{fig:bean_bag_toss} we present the projection of the BRS onto mode 2. The blue dots correspond to the initial values of states on mode 2 that reach $X_T$ at some time for all admissible values of $\mu_k$ (the true BRS); the dashed red line is the outer approximation obtained by utilizing the formulation in $\S$\ref{ssec:free_time}.
\subsection{An example of inner approximation}
Consider the 1D system whose dynamics is described by a variant of the logistic resource growth equation
\begin{align}
  \dot x =\, 0.2x^2+\theta x,\phantom{8} \forall (x,\theta)\in [-1,1]\times [0,0.3].
  \label{eq:ex:1d_inner}
\end{align}
This system exhibits transcritical bifurcation as the value of $\theta$ changes. That is, when $\theta=0$, the origin is attractive for all $x< 0$, and for other values of $\theta$, the origin is unstable and a new stable equilibrium is created in the left half plane; regardless of the value of $\theta$, the origin is repulsive for all $x>0$.
\par
\begin{figure}[!t]
\centering
  {\includegraphics[trim =1.5in 5.5in 1.5in 3.4in, clip=true,width=0.7\columnwidth]{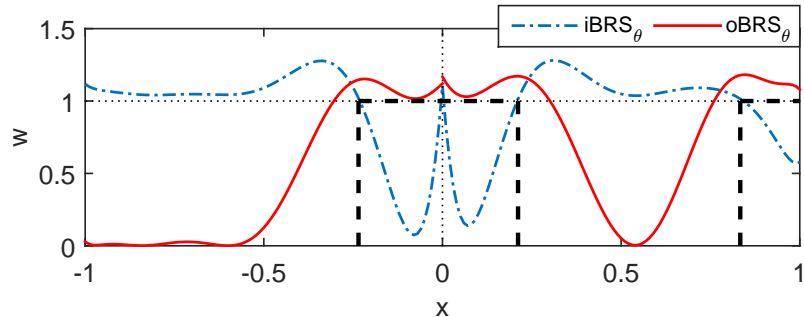}}
  \caption{Outer and inner approximations of the BRS of a 1D system exhibiting transcritical bifurcation. The terminal set, $X_T, T=1$ is a ball of radius 0.3 about the origin. The dashed black boxes are characteristic functions of the true BRS of the uncertain system (Eqn.~(\ref{eq:ex:1d_inner})). The \emph{indicator} functions of degree 12 for outer and inner approximations are presented in solid red and dot-dashed blue respectively.}
  \label{fig:1D_inner}
\end{figure}

Let the space be divided two modes $D_1=[-1,0]$ and $D_2=[0,1]$, with a guard, $G_{(1,2)}$ at $x=1$, and associated reset map $R_{(1,2)}(x)= -x/6$. This reset maps traps trajectories leaving the space from the right and maps them into the left half plane, and closer to the equilibria. Suppose it is of interest to determine the BRS associated with the terminal set $X_T=[-0.3,0.3]$ with $T=1$. Figure~\ref{fig:1D_inner} presents the degree 12 outer (oBRS$_\theta$) and inner approximations (iBRS$_\theta$) of the \emph{indicator function}, $w$, on the BRS, and the true BRS of this uncertain system. Observe that the inner and outer approximations of the BRS are as expected.

%% file: sections/sec_conc.tex
\section{Conclusions}
\label{sec:conclusion}
In this paper, a convex optimization approach is presented to compute the backwards reachable set of quasi-uncertain hybrid systems.
The presented method optimizes over the set of unsigned measures using converging moment relaxations that can be solved using SDPs.
A commentary on the accuracy and the adequacy of the proposed method is provided using examples.
Future work will extend the work herein by incorporating control laws while computing the backwards reachable set.

%% file: sections/sec_app.tex
\section{Existence of solutions}
In Lemma~\ref{lemma:existence} we concern ourselves with establishing that if measures $\mu_{s_{\mathcal J}}$, $\mu_{\mathcal J}$, and $\mu_{f_\mathcal J}$ satisfy Eqn.~(\ref{eq:liouville_1}), then there exist solutions to the hybrid system $\mathcal H$ that originate in the support of $\mu_{s_\mathcal J}$ and terminate at $t=T$ in the support of $\mu_{T_\mathcal J}$. 
\par
The critical ingredient of the proof is the relation between the existence of solutions to a Conservative Continuity Equations (CCE) and an Ordinary Differential Equation (ODE) -- if a solution to the CCE exists, solutions to the ODE exist \cite{Ambrosio2008}. As a review, a continuity equation is a PDE of the following form
\begin{align}
  \frac{\partial \rho}{\partial t}+\nabla\cdot(\rho \mathbf{v})=0,
\end{align}i
where $\mathbf v$ is the `flux' and $\rho$ is the conserved quantity \emph{mass}.
\par
We first show that solutions to CCEs are related to the \emph{average occupation measure} defined in Eqn.~(\ref{eq:mu_avg}). Subsequently, we show that the Liouville eqn. in Eqn.~(\ref{eq:liouville_1}) is a CCE with the associated ODE given by Eqn.~(\ref{eq:ode}), and conclude that $\mu_{\mathcal J}$ which solves the Liouville eqn. is an \emph{average occupation measure} associated with $\mathcal H$. Since $\mu_{\mathcal J}$ solves a CCE, solutions to the ODE in Eqn.~(\ref{eq:ode}) and hence to $\mathcal H$, exist.
\par
As as review, we present, below, the problem of finding solutions to CCE and an associated theorem about the existence of a representing measure.

\begin{lem}
\label{claim:occ}
Consider the homogeneous conservative continuity problem in any one mode of the system; say mode $j$:
\begin{flalign}
&&  \frac{d}{dt}\zeta_t +\nabla\cdot(\tilde f_j\,\zeta_t)=0&&
\label{eq:continuity}
\end{flalign}
where $\zeta_t$ are stochastic kernels supported on $D_j$, conditioned on time, $t$; and $\tilde f_j$ is the Lipschitz vector field. If $\zeta_t$ satisfies the continuity equation, then there exists $\eta\in \mathcal M_+(D_j\times \Xi)$ such that $\forall (\phi,t)\in C_b(D_j)\times \mathcal T$
\begin{align}
  \int_{D_j} \phi(x)\,d\zeta_t=\int_{D_j\times \Xi}\phi(\gamma(t))\,d\eta(x,\gamma),
\end{align}
where $\Xi$ is the space of all absolutely continuous functions supported on $\{\mathcal T\}$ and map to $D_j$; i.e. $\Xi:=C(\mathcal T;D_j)$; and $C_b$ is the space of bounded functions. The measure $\eta$ can be interpreted as a measure on the space of absolutely continuous solutions to the differential equation
\begin{subequations}
\begin{flalign}
  &&\dot\gamma(t) =&\, \tilde f(\gamma(t))\\
  &&\gamma(0) = &\,x_0,& x_0\in D_j
\end{flalign}
\label{eq:app:ODE}
\end{subequations}
such that the conditionals $\eta_x$ are dirac masses.
\par
In addition, if $\lambda$ is the Lebesgue measure on $\mathcal T$, the measure $\zeta^\eta_t\otimes \lambda$ is the average occupation measure defined in Eqn.~(\ref{eq:mu_avg}).
Here $\zeta_t^\eta$ satisfies the following equality
\begin{flalign}
    &&\ip{\zeta_t^\eta,\phi}=\int_{D_j\times \Xi} \phi(\gamma(t))\,d \eta(x,\gamma)&&\forall \phi\in C_b(D_j)
\end{flalign}
\end{lem}
\begin{pf}
The first part of this Lemma is an immediate consequence of Theorem 3.1 in \cite{Ambrosio2008}; the remainder of this proof will establish the relation between the \emph{average occupation measure} and solutions to the continuity equation.\par
  By definition, the following equality holds
  \begin{align}
   (\zeta^\eta_t\otimes\lambda)(A\times B)=&\,\int_{\mathcal T}\int_{D_j\times \Xi}I_{A\times B}(t,x(t))\,d(\zeta^\eta_t \otimes\lambda)
   \end{align}
  Since $\eta\in \mathcal M_+(D_j\times \Xi)$, and this is a Polish space, we can decompose $\eta$ to satisfy the following equality \cite[Corr. 10.4.13]{Bogachev2007a}
  \begin{align*}
    d\eta= d\eta_xd\mu_0
  \end{align*}
  where $\mu_0$ is the distribution of initial conditions for the ODE in Eqn.~(\ref{eq:app:ODE}) and $\eta_x$ is a regular conditional measure, conditioned on the initial value of $x$. Thus,
   \begin{align}
       (\zeta^\eta_t\otimes\lambda_t)(A\times B)=&\,\int\limits_{\mathcal T}\hspace*{-2mm}\int\limits_{D_j}\hspace*{-2mm}\int\limits_{\Xi}\hspace*{-2mm}I_{A\times B}(t,x(t\mid x_0))\,d\eta_x\,d\mu_0\,dt\label{eq:app:claim:1}\\
    =&\,\int\limits_{\mathcal T}I_A(t)\int\limits_{D_j}I_B(x(t\mid x_0))\,d\mu_0\,dt\label{eq:app:claim:2}\\
    =&\,\int\limits_{D_j}\int\limits_{\mathcal T}I_{A\times B}(t,x(t\mid x_0))\,dt\,d\mu_0\label{eq:app:claim:3}\\
    =&\,\int\limits_{D_j}\mu_j(A\times B\mid x_0)\,d\mu_0\label{eq:app:claim:4}\\
    =&\, \mu_j(A\times B)\label{eq:app:claim:5};
  \end{align}
where the transition from Eqn.~(\ref{eq:app:claim:1}) to Eqn.~(\ref{eq:app:claim:2}) uses the fact that $\eta_x$ is a dirac measure; Fubini's theorem is used to get to Eqn.~(\ref{eq:app:claim:3}); and the definition of \emph{occupation measure} and  \emph{average occupation measure} are used to arrive at Eqns.~(\ref{eq:app:claim:4}) and (\ref{eq:app:claim:5}) respectively.
\end{pf}

\begin{lem}
Let $\big(\mu_{0_{\mathcal J}},\mu_{T_{\mathcal J}},\mu_{G_\mathcal E},\mu_{\mathcal J}\big)$ satisfy Eqn.~(\ref{eq:primal:liouville}). Then, there exists a family of absolutely continuous state trajectories starting emanating from each $\mu_{0_j},\,j\in \mathcal J$ such the occupation and terminal measures (at $t=T$) in each mode generated by this family of trajectories coincide with $\mu_j$ and $\mu_{T_j},\,j\in \mathcal J$. In addition, $\mu_{G_\mathcal E}$ s coincide with the restrictions of $\mu_{\mathcal J}$ on the respective guards.
\label{lemma:existence}
\end{lem}
\begin{pf}
This lemma is a generalization of \cite[Lemma 3]{henrion2014convex} and is proven using the same technique adopted by the authors of \cite{henrion2014convex}. In presenting this proof, it is assumed that uncertain parameters are augmented to the state-space description of the system, similar to $\S$\ref{sec:prob}. In addition, given the description of hybrid systems considered in this paper (refer $\S$\ref{ssec:hybrid}), it is sufficient to restrict our attention to one particular mode of the system, say mode $j$. For notational convenience, let $\sigma_{0_j}$ denote the measure on the image of reset maps originating from the guards of in $\mathcal H$; and let the starting measure $\mu_{s_j}$ be defined as
\begin{align}
\mu_{s_j} = \delta_0\otimes\mu_{0_j}\otimes\mu_{\theta_j}+\sigma_{0_j}
\label{eq:app:musj}
\end{align}
In accordance with the hybrid system definition, trajectories of the state can arrive in mode $j$ at time $\tau\in [0,T]$ via $\spt(\mu_{s_j})$; and once in the mode, can arrive at $\spt(\mu_{f_j})$ at any time $t\in [0,T]$, and hence possibly leave the mode. Equation~(\ref{eq:liouville_1}) has a conserved quantity -- \emph{mass}. To see this, consider the test function $\phi(t) = t^k$ with $k=0$; then it follows that
\begin{flalign}
  \mu_{s_j}(A\times B\times C) = \mu_{f_j}(A\times B\times C).
  \label{eq:app:def:sj}
\end{flalign}
That is, along the flow of solutions in each mode, \emph{mass} is conserved and the Liouville eqn. is a CCE in terms of the measure $\mu_j$. The measures given by the statement of this Lemma satisfy the Liouville eqn. and hence a CCE; thus solutions to the associated family of ODEs exists in each mode. To see that solutions to the hybrid description of ODEs exists, refer to the results in \cite{Burden2015}. Also, using Lemma~\ref{claim:occ}, since $\mu_j(x,\theta\mid t)$ solves the CCE in each mode, $\mu_j$ is the \emph{average occupation measure}.
\par
Having established that solutions to the ODE exists in every mode, we have to demonstrate that the solution trajectories begin in $\bigcup_{j\in \mathcal J}\spt(\mu_{0_j})$ and terminate in $\bigcup_{j\in \mathcal J}\spt(\mu_{T_j})$; this, again, can be demonstrated in a per-mode basis. Since $\mu_j(\cdot\mid t)$ solves a CCE, conditioned at $t=0$, this measure coincides with the initial distribution ($\nu_{0_j}$) of states and at $t=T$, it coincides with the final distribution of states; we have to show that $\spt(\nu_{0_j})=\spt(\mu_{0_j})$ and that $\spt(\mu_j(\cdot\mid T))=X_{(T,j)}$. We do the same in the ensuing presentation.
\par
Since we are concerned with solutions evolving on Polish spaces, we can decompose $\mu_j$ into the following form \cite[Corr. 10.4.13]{Bogachev2007a}
\begin{align}
  d\mu_{j}(t,x,\theta)=d\tilde\mu_{j}(x,\theta\mid t)\, d\xi_{\mu_j}(t)
  \label{eq:app:mu:decomp}
\end{align}
where $\tilde \mu_j(\cdot\mid t)$ is the regular conditional measure and $\xi_{\mu_j}$ is the normalized projection of $\mu_j$ on $\mathcal T$. We claim that if $\lambda$ is the Lebesgue measure on $\mathcal T$, $\xi_{\mu_j}$ is a scaled version of $\lambda$. That is, we have to demonstrate that the moments of $\xi_{\mu}$ are a scaled version of that of $\lambda$. Since each $\mu_j$ is an \emph{average occupation measure}, by definition (Eqn.~(\ref{eq:mu_avg})), we have that
\begin{align}
\mu_j(A\times B\times C)=&\,\hspace*{-5mm}\int\limits_{\mathcal T\times D_j\times \Theta_j}\hspace*{-6mm}\int_{0}^T \hspace*{-2.25mm} I_{A\times B\times C}(t,x(t|\tau_k,x_0,\theta),\theta)\,dt\,d\nu_{0_j},
\label{eq:app:mu:lebesgue}
\end{align}
where $\nu_{0_j}\in \mathcal M_+(\mathcal T\times D_j\times \Theta_j)$ is the initial distribution of states (not yet known to be related to $\mu_{s_j}$). Hence, the t-moments of $\mu_j$ are
\begin{align}
 \int\limits_{\mathcal T\times D_j\times \Theta_j}t^k \,d\mu_j=&\,\int\limits_{\mathcal T\times D_j\times \Theta_j} \int\limits_{0}^T \hspace*{-1.25mm} t^k \,dt\,d\nu_{0_j},\\
 =&\,\frac{T^{k+1}}{k+1}\nu_{0_j}(\mathcal T\times D_j\times \Theta_j),
\end{align}
which are scaled moments of the $\lambda$. Thus it follows that Eqn.~(\ref{eq:app:mu:decomp}) can be written as
\begin{align}
  d\mu_{j}=d\mu_{j}(x,\theta\mid t) \,dt.
  \label{eq:app:mu:decomp:final}
\end{align}
Recall that in each mode, the supports of measures in Liouville eqn. have the following properties: $\spt(\mu_{s_j})\subset (\mathcal T\times D_j^\circ\times \Theta_j)$; $\spt(\mu_j)\subset (\mathcal T\times D_j\times \Theta_j)$, $\spt(\mu_{T_j})\subset X_{(T,j)}\backslash \bigcup_{e:=(j,k),e\in \mathcal E} G_{e}\times \Theta_j$; and $\spt(\mu_{G_e})\subset (\mathcal T\times G_e\times \Theta_j)$. Observe that only $\mu_j$ and $\mu_{G_e}$ are supported on $\mathcal T\times G_e\times \Theta_j$; it thus follows from Eqn.~(\ref{eq:primal:liouville})that for all $e:=(j,k),\forall k\in \{l\mid (j,l)\in \mathcal E\}$
\begin{align}
  \mu_{G_e}(\mathcal T\times G_e\times\Theta_j) = \mathcal L_{f_j}'\mu_j(\mathcal T\times G_e\times \Theta_j).
  \label{eq:app:mu_G_arg1}
\end{align}
Since $\mu_j$ is an unsigned measure whose t-marginal is a scaled version of $\lambda$, it follows that for all $A\subset \mathcal T\text{ st. }\lambda(A)=0$
\begin{flalign}
  &&0=&\,\mu_j(A\times D_j\times \Theta_j)\ge \mu_j(A\times G_e\times \Theta_j).&
  \label{eq:app:mu_G_arg2}
\end{flalign}
From Eqns.~(\ref{eq:app:mu_G_arg1}) and (\ref{eq:app:mu_G_arg2}), it follows that the t-marginal of $\mu_{G_e}$ is absolutely continuous wrt. $\lambda$. Thus, the measures on the guards can be decomposed as follows:
\begin{flalign}
  &&d\mu_{G_{e}}=&\,w_{e}(t)\,d\mu_{G_e}(x,\theta\mid t)\, dt, & \forall e\in \mathcal E
  \label{eq:app:mu_G_arg3}
\end{flalign}
where $w_e(t)$ is the density of the t-marginal of $\mu_{G_e}$ wrt. to $\lambda$.
\newline
In each mode $j\in \mathcal J$, $\sigma_{0_j}$ is a sum of the push-forward of some measures on guards through reset maps that are identity maps in the t-component; that is,
\begin{align}
\sigma_{0_j}=\sum_{e\in \{(k,j)\mid \mathcal E\}}R^*_{e}(\pi^*_{(t,x)}\mu_{G_e})\otimes \mu_{\theta_j}.
\end{align}
Now, by considering test functions of the form $v(t,x,\theta)=\varphi(t)\phi(x,\theta)$ in Eqn.~(\ref{eq:primal:liouville}), we get (by employing integration by parts):
  \begin{align}
    \varphi(T)\alpha_T-\varphi(0)\alpha_0=\int_0^{T} \dot \varphi \alpha_t+\varphi\cdot(\alpha_f+\alpha_\sigma-\alpha_G)\,dt
        \label{eq:app:2}
  \end{align}
  where
  \begin{flalign}
      &&&\alpha_T:=\,\int\limits_{X_{(T,j)}\times \Theta_j}\hspace*{-.1in}\left[\phi(x,\theta)\right]\,d\mu_{T_j},&\\
      &&&\alpha_0:=\,\int\limits_{X_{(0,j)}\times \Theta_j}\hspace*{-.1in}\left[\phi(x,\theta)\right]\,d \mu_{0_j}d\mu_{\theta_j}, \\
      &&&\alpha_t(t):=\,\int\limits_{{D}_j\times \Theta_j}\left[\phi(x,\theta)\right]\,d\mu_{j}(x,\theta\mid t)\label{eq:app:c},\\
      &&&\alpha_f(t):=\,\int\limits_{{D}_j\times \Theta_j}\left[\nabla\phi(x,\theta)\cdot f_j\right]\,d\mu_{j}(x,\theta\mid t), \\
      &&&\alpha_\sigma(t):=\,\hspace*{-5mm}\sum_{e\in \{(k,j)\in \mathcal E\}}\hspace*{.25mm}\int\limits_{G_e\times\Theta_k}\hspace*{-3mm}\left[w_e\int_{\Theta_j}\phi(R_e(x),\theta)\,d\mu_{\theta_j}\right]\,d\mu_{G_e}(x,\theta\mid t),\\
      &&&\alpha_{G}(t):=\,\sum_{e\in \{(j,k)\in \mathcal E\}}\int_{G_e\times \Theta_j}\left[w_e \phi(x,\theta)\right]\,d\mu_{G_e}(x,\theta\mid t).
      \end{flalign}
Given functions $\alpha_T,\alpha_0,\alpha_f,\alpha_\sigma$ and $\alpha_G$, Eqn.~(\ref{eq:app:2}) admits an $dt$-a.e. unique solution for $\alpha_t$. To see this, suppose $c_1(t)$ and $c_2(t)$ are any two admissible solutions. Then
\begin{subequations}
\begin{align}
    \varphi(T)\alpha_T-\varphi(0)\alpha_0=&\,\hspace*{-1mm}\int\limits_0^{T} \hspace*{-1mm}\left[\dot \varphi  c_1+\varphi(\alpha_f+\alpha_\sigma-\alpha_G)\right]dt\label{eq:app:2.5a}\\
    \varphi(T)\alpha_T-\varphi(0)\alpha_0=&\,\hspace*{-1mm}\int\limits_0^{T} \hspace*{-1mm}\left[\dot \varphi  c_2+\varphi (\alpha_f+\alpha_\sigma-\alpha_G)\right]\,dt\label{eq:app:2.5b}
\end{align}
\end{subequations}
Taking the difference between Eqns.~(\ref{eq:app:2.5a}) and (\ref{eq:app:2.5b}) and noting that $\varphi(t)\in \mathcal C(\mathcal T)$ is any arbitrary function and that $\mathcal C(\mathcal T)$ is dense in $L^1(\mathcal T)$, it follows that $c_1(t)=c_2(t)$, $dt$-a.e..
\par
By construction, with $\alpha_T,\alpha_0,\alpha_f,\alpha_\sigma$ and $\alpha_G$ as defined, we know that $\alpha_t(t)$ as defined in Eqn.~(\ref{eq:app:c}) satisfies Eqn.~(\ref{eq:app:2}). We claim that $\tilde \alpha_t(t)$ defined as
\begin{align}
    \tilde \alpha_t(t) = \alpha_0+\int\limits_{0}^t \left[\alpha_f+\alpha_\sigma-\alpha_G\right]\,d\tau
    \label{eq:app:3}
\end{align}
also solves Eqn.~(\ref{eq:app:2}). Observe from Eqn.~(\ref{eq:app:3}) that:
\begin{align}
\tilde \alpha_t(0)   = &\, \alpha_0+ \int\limits_{0}^0 \left[\alpha_f+\alpha_\sigma-\alpha_G\right]\,d\tau = \alpha_0,\\
\tilde \alpha_t(T)   = &\, \alpha_0+ \int\limits_{0}^T \left[\alpha_f+\alpha_\sigma-\alpha_G\right]\,d\tau = \alpha_T\label{eq:app:alphaT}\\
\frac{d\tilde \alpha_t}{dt} = &\, \alpha_f+\alpha_\sigma-\alpha_G,
\end{align}
where the equality in Eqn.~(\ref{eq:app:alphaT}) follows from considering $v=\varphi(x,\theta)$ as the test function in Eqn.~(\ref{eq:primal:liouville}), and the last equality follows from Leibniz rule.
\par
Now consider evaluating $\int_{\mathcal T}(\dot\varphi\tilde \alpha_t)\,dt$.
\begin{align}
  \int\limits_{\mathcal T}(\dot\varphi\tilde \alpha_t)\,dt = \varphi(T)\tilde \alpha_t(T)-\varphi(0)\tilde \alpha_t(0)-\int\limits_0^{T} \varphi\frac{d\tilde \alpha_t}{dt}\,dt\\
  = \varphi(T)\alpha_T-\varphi(0)\alpha_0 -\int\limits_0^{T}\varphi\left[\alpha_f+\alpha_\sigma-\alpha_G\right]\,dt.
\end{align}
This is of the form of Eqn.~(\ref{eq:app:2}) and hence $\tilde \alpha_t(t)$ solves Eqn.~(\ref{eq:app:2}).
\par
Since the solution to Eqn.~(\ref{eq:app:2}) is $dt$-a.e. unique, it follows that, $\forall \phi\in \mathcal C^1(D_j\times \Theta_j)$:
\begin{align}
    \int\limits_{D_j\times \Theta_j}\phi(x,\theta)\,d\mu_j(x,\theta\mid t)=&\,\alpha_0+\int\limits_{0}^t \left[\alpha_f+\alpha_\sigma-\alpha_G\right]\,d\tau
    \label{eq:app:unique_sol}
\end{align}
From Eqn.~(\ref{eq:app:unique_sol}) it follows that the stochastic kernel of $\mu_j$ at $t=0$ coincides with $\mu_{0_j}$, and that it coincides with $\mu_{T_j}$ at $t=T$. This implies that the solutions to the associated ODE, at time $t=T$ terminate in $\spt(\mu_{T_j})=X_{(T,j)}$.
\par
To complete the proof, we now show that the restriction of $\mu_j$ onto the guards is, $dt$-a.e., identical to the measure on the guards; further that $\nu_0$ (from Eqn.~(\ref{eq:app:mu:lebesgue})) is identical to $\mu_{s_j}$ (from Eqn.~(\ref{eq:app:def:sj})).
\par
Recall the definition of the average occupation measure (Eqn.~(\ref{eq:mu_avg})), its relation to the initial state distribution $\mu_{s_j}$, and the form of the Liouville eqn. (Eqn.~(\ref{eq:primal:liouville})); it is clear that (in mode $j$) the following equality holds
\begin{align}
  \nu_0+\mathcal L_{f_j}'\mu_j=\delta_T\otimes\mu_T+\mu_G,
  \label{eq:app:muG:arg1}
\end{align}
where $\mu_G$ is the sum of measures on all guards in mode $j$.\par
From Remark~\ref{rem:mu_G:restriction}, the following equality holds
\begin{align}
  \nu_0+\mathcal L_{f_j}'\mu_j=\delta_T\otimes\mu_T+\sum_{k\in \{l\mid (j,l)\in \mathcal E\}}\mu_{j}|_{G_{(j,k)}}.
  \label{eq:app:muG:arg2}
\end{align}
Taking the difference between Eqns.~(\ref{eq:app:muG:arg1}) and (\ref{eq:app:muG:arg2}), noting that since $\mu_j$ is decomposable in the form of Eqn.~(\ref{eq:app:mu:decomp:final}), so is its restriction, recalling that $\mu_G$ is also decomposable (Eqn.~(\ref{eq:app:mu_G_arg3})), and considering the product of arbitrary test functions in $\mathcal C^1(\mathcal T)$ and $\mathcal C^1(D_j\times \Theta_j)$, it is concluded that the t-conditionals of $\mu_G$ and $\mu_j$ are equal $dt$-a.e.. Hence, $\mu_G$ is the restriction of $\mu_j$ to the guards. That that $\nu_0$ is equal to $\mu_{0_j}+\mu_{s_j}$ follows as a natural consequence of Eqn.~(\ref{eq:app:musj}) and our conclusion about the t-conditional of $\mu_j$ at $t=0$.
\end{pf}